\documentclass[12pt]{article}
\usepackage{colonequals}
\usepackage{booktabs}
\usepackage{subfigure}
\usepackage{enumerate}
\usepackage{amsmath}
\usepackage{amssymb}
\usepackage{amsfonts}
\usepackage{multirow}
\usepackage{amsthm}
\usepackage{cmll}
\usepackage{authblk}
\usepackage{geometry}
\usepackage{natbib}
\usepackage{array,epsfig,fancyheadings,rotating}
\newcommand{\var}{{\rm var}}
\newtheorem{theorem}{\bf{Theorem}}

\newtheorem{example}{\bf{Example}}

\newtheorem{condition}{\bf{Condition}}

\newtheorem{assumption}{\bf{Assumption}}

\def\T{{ \mathrm{\scriptscriptstyle T} }}

\begin{document}

	\def\spacingset#1{\renewcommand{\baselinestretch}%
		{#1}\small\normalsize} \spacingset{1}

	%%%%%%%%%%%%%%%%%%%%%%%%%%%%%%%%%%%%%%%%%%%%%%%%%%%%%%%%%%%%%%%%%%%%%%%%%%%%%%
	
	\title{Sharp  bounds for  variance of treatment effect estimators in the finite population in the presence of covariates}
	%% The year, volume, and number are determined on publication
	\author{Ruoyu Wang\thanks{wangruoyu17@mails.ucas.edu.cn}} 
	\author{Qihua Wang
		\thanks{qhwang@amss.ac.cn}}
	\affil{Academy of Mathematics and Systems Science, Chinese Academy of Sciences, 55 Zhongguancun East Road, Haidian District, Beijing 100190, China}
	
	\author{Wang Miao
		\thanks{mwfy@pku.edu.cn}}
	\affil{School of Mathematical Sciences, Peking University, 5 Summer Palace Road, Haidian District, Beijing 100871, China}
	
	\author{Xiaohua Zhou
		\thanks{azhou@bicmr.pku.edu.cn}}
	\affil{Beijing International Center for Mathematical Research and Department of Biostatistics, Peking University, 5 Summer Palace Road, Haidian District, Beijing 100871, China
	}
	\maketitle

\begin{abstract}
	In a completely randomized experiment, the variances of treatment effect estimators in the finite population are usually not identifiable and hence not estimable.   Although some estimable bounds of the variances have been established in the literature,  few of them are derived in the presence of covariates. 
	In this paper, the difference-in-means estimator and the Wald estimator are considered in the completely randomized experiment with perfect compliance and noncompliance, respectively.  Sharp bounds for the variances of these two estimators are established when covariates are available.
	Furthermore,   consistent estimators for such bounds are obtained, which can be used to shorten the confidence intervals and improve the power of tests. Confidence intervals are constructed based on the consistent estimators of the upper bounds, whose coverage rates are uniformly asymptotically guaranteed.
	Simulations were conducted to evaluate the proposed methods. The proposed methods are also illustrated  with two real data analyses.
	
\end{abstract}
\noindent%
\vfill

\section{Introduction}
Estimation and inference for the average treatment effect are extremely important in practice. 
Lots of the literature assumes that the observations are sampled from an infinite super population \citep{hirano2003efficient,Imbens2004,belloni2014inference,chan2016globally}. The infinite super population seems contrived if we are interested in evaluating the treatment effect for a particular finite population \citep{Li2017GFoFPCLTwAtCI}, e.g., the patients enrolled in the experiment,  and the finite population framework is more suitable for such problems. 
The finite population framework views all potential outcomes as fixed and the randomness of the data comes from the treatment assignment solely \citep{Imbens2005,Nolen2011RBIwPS}. This framework avoids assumptions about randomly sampling from some ``vaguely defined super-population of study units” \citep{Schochet2013cluster}.  The statistical analysis results under the finite population framework are interpretable in the absence of any imaginary super-population. Theoretical guarantees in this framework rely on the treatment assignment rather than unverifiable sampling assumptions such as the independent and identically distributed assumption. In randomized experiments, 
%which are deemed as a gold standard for evaluating treatment effects, 
the finite population framework has been widely used in data analysis since \citet{neyman1990application}. A fundamental problem in the completely randomized experiment under the finite population framework is that the variance of the widely-used difference-in-means estimator is unidentifiable.
Thus, we can not obtain a consistent variance estimator and the standard inference based on normal approximation fails.
To mitigate this problem,  \citet{neyman1990application} 
%decomposed the variance and omitted a negative term to derive an 
adopted an estimable upper bound for the variance, which leads to conservative inference. 
%Replacing the variance with a consistent estimation of its upper bound in the inference procedure leads the inference to be conservative. 
The precision of the bound is crucial for the power of a test and the width of the resulting confidence interval. 
Thus it is important to incorporate all the available information to make the bound as precise as possible. The variance bound with binary outcomes  was fairly well studied by \citet{Robins1988CIfCP}, \citet{Ding2016APToTbTTfCRE},  \citet{Ding2018MFCIoBED}, etc.
For general outcomes, \cite{Aronow2014SBotViRE} improved \citet{neyman1990application}'s results by deriving a  sharp bound that can not be improved without information other than the marginal distributions of potential outcomes.

In many randomized experiments, some covariates are observed in addition to the outcome.
However, few previous approaches consider how to improve the variance bound using covariate information, with the exception of \cite{Ding2018DTEV}. 
Surprisingly, we observe that  the upper bound for variance of the treatment effect estimator given by \cite{Ding2018DTEV}  can be larger than that given by \citet{Aronow2014SBotViRE} in some situations. 
This is illustrated with Example \ref{ex: bound perform} in Section \ref{sec:CRT}. 
%In a completely randomized experiment, according to Neymanian decomposition, the %variance of the difference-in-means estimator is decomposed into two identifiable terms %and an unidentifiable term.

The first main contribution of this paper is to derive the sharp bound for the variance of the difference-in-means estimator in the finite population when covariates are available and obtain  a consistent estimator of the bound. The proof of  
consistency is quite challenging.  In the analysis of consistency, we allow the cardinality of the covariate support to diverge as the population size. This is different from that in the literature and increases the difficulty 
of the proof due to the lack of tools for analyzing sample conditional quantile functions involved in the estimator in finite population when the cardinality of the covariate support is diverging.
Based on the consistent estimator of the variance bound, a shorter confidence interval with a more accurate coverage rate is obtained. In addition, we show that the confidence interval has an asymptotically guaranteed coverage rate and the asymptotic result is uniform over a large class of finite populations. As discussed by \citet{lehmann2006testing}, the uniformity is crucial for reassuring inference based on asymptotic results and has long been omitted by existing works.

The aforementioned results focus on completely randomized experiments where units comply with the assigned treatments.
%with perfect compliance where the treatment a unit is assigned is the treatment the unit actually receives. 
However, noncompliance often occurs in randomized experiments, 
in which case, the parameter of interest is the local average treatment effect  \citep{angrist1996identification,abadie2003semiparametric}, 
and its inference is more complicated.
%In the presence of  noncompliance, the problem becomes more complicated. 
%In this case, rather than the average treatment effect, the parameter of interest is the local average treatment effect
% under different randomization schemes (both are different from the randomization scheme we consider). 
Recently,  the estimating problem of the local average treatment effect is considered in the finite population. Some estimators are suggested including the widely-used Wald estimator and the ones due to  \cite{Ding2018DTEV} and \cite{Hong2020}. Identification problem also exists for the variances of these estimators.  
However, to the best of our knowledge, no sharp bound is obtained for the unidentifiable variances in the literature.  

Another main contribution of this paper is to further  extend the aforementioned results for the
completely randomized experiment to the case 
with noncompliance. We establish the sharp bound for the variance of the Wald estimator and propose the consistent estimator for the variance bound. The analysis of consistency is even more involved in this case due to the complexity of the estimator. A confidence interval based on the consistent estimator of the upper bound is also constructed, whose coverage rate is uniformly asymptotically guaranteed. It is worth mentioning that the sharp bound without covariates can be derived as a special case of the resulting bound, which has not been investigated in the literature. 
%Construction of the bounds and proofs of consistency are novel in the finite population framework. 
Simulations and application to two real data sets from the randomized trial ACTG protocol 175 \citep{hammer1996} and JOBS II \citep{vinokur1995} demonstrate the advantages of our methods.

This paper is organized as follows.  In Section \ref{sec:CRT}, we establish the sharp variance bound in the presence of covariates for the difference-in-means estimator in the completely randomized experiment with perfect compliance. A consistent estimator is obtained for the bound. In Section \ref{sec:noncom}, we consider the Wald estimator for the local average treatment effect in the completely randomized experiment in the presence of noncompliance; we establish a sharp variance bound for the Wald estimator in the presence of covariates and obtain a consistent estimator for the bound. Simulation studies are conducted to evaluate the empirical performance of the proposed bound estimators in Section \ref{sec:simulation}, followed by some applications to data from the randomized trial ACTG protocol 175 and JOBS II in Section \ref{sec:real data}. A discussion on some possible extensions of our results is provided in Section \ref{sec:discuss}. Proofs are relegated to the Appendix.

%\section{Preliminaries}\label{subsec:preliminaries}

%In what follows, we derive the sharp bound and for the variance of treatment effect 
%Throughout this paper, we make the convention $0/0=0$.

\section{Sharp variance bound for the difference-in-means estimator }\label{sec:CRT}
\subsection{Preliminaries}\label{subsec:preliminaries}

Suppose we are interested in the effect of a binary treatment on an outcome in a finite population consisting of $N$ units.
In a completely randomized experiment, $n$  out of $N$ units are sampled from the population, 
with $n_1$ of them randomly assigned to the treatment group and  the other  $n_0=n-n_1$ to the control group. 
Let $T_i=1$ if unit $i$ is assigned to the treatment group, $T_i=0$ the control group, and $T_i$ is not defined if unit $i$ is not enrolled in the experiment. 
For each unit $i$ and $t=0,1$,  let $y_{ti}$ denote the  potential outcome that would be observed if unit $i$ is assigned to treatment $t$.
Let   $w_i$  denote a vector  of covariates with the constant $1$ as its first component. The covariate vector $w_{i}$ is observed if the unit $i$ is enrolled in the experiment (i.e., $T_{i} = 0$ or $1$).
Then the characteristics of the population  can be viewed as a matrix $\mathbf{U}=(y_1,y_0,w)$ where $y_1=(y_{11},y_{12},\dots,y_{1N})^T$, $y_0=(y_{01},y_{02},\dots,y_{0N})^T$
and $w=(w_1,\dots,w_N)^\T$. 

For any vector $a=(a_1,\dots,a_N)^\T$, we let
\[\mu(a)=\frac{1}{N}\sum\limits_{i=1}^Na_i,\quad \phi^2(a)=\frac{1}{N}\sum\limits_{i=1}^N(a_i-\mu(a))^2.\]
Letting $\tau_i=y_{1i}-y_{0i}$ be the treatment effect for unit $i$ and $\tau=(\tau_{1},\dots,\tau_N)^\T$,  the parameter of interest is the  average treatment effect, 
\[\theta = \mu(\tau) = \frac{1}{N}\sum\limits_{i=1}^Ny_{1i} - \frac{1}{N}\sum\limits_{i=1}^Ny_{0i}.\] 
Note that all the parameters in this paper depend on $N$ if not otherwise specified, and we drop out  the dependence in notation for simplicity when there is no ambiguity.
The treatment assignment is unrelated to the covariates in completely randomized experiments. Hence the average treatment effect can be estimated by the difference-in-means estimator
\[\hat{\theta}=\frac{1}{n_1}\sum_{T_i=1} y_{1i}  - \frac{1}{n_0}\sum_{T_i=0} y_{0i},  \]
This estimator is widely used due to its simplicity, transparency, among other practical reasons \citep{shao2010theory,Lin2013ANoRAtEDRFC}. Moreover, it is the uniformly minimum variance unbiased estimator under the scenario presented in \citep{kallus2018optimal}. As in a lot of the literature \citep{imai2008variance,Aronow2013ACoUEotATEiRE,shao2010theory,kallus2018optimal,ma2020regression}, we consider the inference based on the difference-in-means estimator because of its popularity in practice and its theoretical importance.

According to \cite{Freedman2008ORAiEwST},  the variance of $\hat{\theta}$ is 
\[ \frac{1}{N-1}\left\{\frac{N}{n_1}\phi^2(y_1)+\frac{N}{n_0}\phi^2(y_0)-\phi^2(\tau)\right\},\]
and we denote this variance by $\sigma^{2}/(N-1)$.
Under certain standard regularity conditions in the finite population, previous authors \citep{Freedman2008ORAiEwST,Aronow2014SBotViRE,Li2017GFoFPCLTwAtCI} have established that
\begin{equation}\label{eq: asy-crt}
	\sqrt{N}\sigma^{-1}(\hat{\theta} - \theta) \overset{d}{\to} N(0, 1),
\end{equation}
as $n_{1}$, $n_{0}$ and $N$ goes to infinity.
Statistical inference may be made based on this asymptotic distribution.
However,  it is difficult  to obtain  a consistent  estimator for $\sigma^{2}$. 
According to standard results in survey sampling \citep{cochran1977}, $\phi^2(y_t)$  can be consistently estimated by
\begin{equation}\label{eq: est phi_hat_t}
	\hat{\phi}_{t}^{2} = \frac{1}{n_t - 1}\sum_{T_i=t}\Big(y_{ti} - \frac{1}{n_t}\sum_{T_j=t}y_{tj}\Big)^2
\end{equation}
for $t=0,1$.
However, $\phi^2(\tau)$ and hence $\sigma^2$ is not identifiable 
because the potential outcomes $y_1$ and $y_0$ can never be observed simultaneously.
To make inference for $\theta$ based on (\ref{eq: asy-crt}), one can   use an  upper bound for $\sigma^{2}$  to construct a conservative confidence interval.
Alternatively, one may use an estimable lower bound for $\sigma^2$ to get a shorter confidence interval. However, the coverage rate of such a confidence interval may not be guaranteed. 
To establish an estimable upper (lower) bound for $\sigma^2$, it suffices to establish an estimable lower (upper) bound for 
the unidentifiable term $\phi^{2}(\tau)$. 
We then derive the sharp bound for $\phi^2(\tau)$ and obtain its consistent estimator.

\subsection{Sharp bound for $\phi^2(\tau)$}
For any matrices $a=(a_1,\dots,a_N)^\T$, $b=(b_1,\dots,b_N)^\T$ and vectors $\bar{a}$, $\bar{b}$ whose dimensions equal to the number of columns of $a$ and $b$, respectively, define 
\begin{align*}
	&P(a \leq \bar{a}) = \frac{1}{N}\sum_{i=1}^{N}1\{a_{i} \leq \bar{a}\},\\
	&P(a = \bar{a}) = \frac{1}{N}\sum_{i=1}^{N}1\{a_i = \bar{a}\}, 
\end{align*}
\begin{align*}
	&P(a= \bar{a}\mid b=\bar{b})=\frac{\sum_{i=1}^{N}1\{a_i = \bar{a},b_i =  \bar{b}\}}{\sum_{i=1}^{N}1\{b_i = \bar{b}\}}, \\
	&P(a\leq \bar{a}\mid b=\bar{b})=\frac{\sum_{i=1}^{N}1\{a_i\leq \bar{a},b_i =  \bar{b}\}}{\sum_{i=1}^{N}1\{b_i = \bar{b}\}},
\end{align*}
where $1\{\cdot\}$ is the indicator function and the ``$\leq$" between two vectors corresponds to the component-wise inequality.
Note that in this paper, $P(\cdot)$ and $P(\cdot\mid \cdot)$ are some quantities that describe a vector, and we use $\mathbb{P}(\cdot)$ to denote the probability.

For any function $H$, we define $H^{-1}(u)=\inf\{s:H(s)\geq u\}$. In this paper, we adopt the convention $\inf \emptyset = \infty$.
%A trivial lower bound for $\phi^{2}(\tau)$ is zero, which does not utilize any information from data.  
We let $\{\xi_1,\dots,\xi_K\}$ be the set of all different values of $w_i$.  Clearly $K\leq N$. We aim to derive bounds for  $\phi^2(\tau)$  by using covariate information efficiently. Define  $\pi_k=P(w=\xi_k)$ and for $t=0,1$ and $k=1,\dots,K$. 
The quantities $\pi_{k}$ ($k=1,\dots, K$) and functions $F_{t\mid k}$ ($t=0,1$ and $k=1,\dots, K$) summarize the characteristics of the population and can be estimated using observed data. To obtain estimable bounds for $\phi^{2}(\tau)$,
we focus our attention on the bound that can be expressed as a functional of $\pi_{k}$ and $F_{t\mid k}$ ($t=0,1$ and $k=1,\dots, K$). Define the set of lower bounds 
\begin{align*}
	\mathcal{B}_{\rm L} = \{b_{\rm L}: &\text{$b_{\rm L}$ is a functional of $\pi_{k}$ and $F_{t\mid k}$ for $t=0,1$ and $k=1,\dots, K$;}\\
	&b_{\rm L} \leq \phi^{2}(\tau)\}.
\end{align*}
Define the set of upper bounds $\mathcal{B}_{\rm H}$ similarly.
%Using $F_{t\mid k}$ and $\pi_k$, $k=1,\dots,K$, we derive a lower (upper) bound for $\phi^2(\tau)$ by the following procedure:
%It can be verified that the data distributions of the treatment and the control group under any population in $\mathcal{U}$ are identical to those under $\mathbf{U}$. Thus we can not distinguish $\mathbf{U}$ from other populations in $\mathcal{U}$ based on the data distributions of the treatment and the control group.
Then the sharp bound is established in the following theorem.
\begin{theorem}\label{thm:cov}
	A bound for  $\phi^2(\tau)$ is $[\phi^2_{\rm L},\phi^2_{\rm H}]$ where
	\begin{align*}
		&\phi^2_{\rm L}=\sum_{k=1}^{K}\pi_k\int_0^{1}(F_{1\mid k}^{-1}(u)-F_{0\mid k}^{-1}(u))^2du-\theta^2,\\
		&\phi^2_{\rm H}=\sum_{k=1}^{K}\pi_k\int_0^{1}(F_{1\mid k}^{-1}(u)-F_{0\mid k}^{-1}(1-u))^2du-\theta^2.
	\end{align*}
	Moreover, the bound is sharp in the sense that $\phi_{\rm L}^{2}$ is the largest lower bound in $\mathcal{B}_{\rm L}$ and  $\phi_{\rm H}^{2}$ is the smallest upper bounds in $\mathcal{B}_{\rm H}$.
\end{theorem}

See the Appendix for the proof of this theorem.
Here we compare this bound to previous bounds obtained by  \cite{Aronow2014SBotViRE} and \cite{Ding2018DTEV}.
By utilizing  the marginal distributions of potential outcomes, \citet{Aronow2014SBotViRE} derived the bound for $\phi^2(\tau)$:
\begin{align*}
	\phi^2_{\rm AL}\colonequals&\int_0^{1}(F_1^{-1}(u)-F_0^{-1}(u))^2du-\theta^2\leq\phi^2(\tau)
	\\&\leq\int_0^{1}(F_1^{-1}(u)-F_0^{-1}(1-u))^2du-\theta^2\equalscolon\phi^2_{\rm AH}
\end{align*}
where $F_t(y)=P(y_t\leq y)$ for $t=0,1$.
The bound of \cite{Aronow2014SBotViRE}  is sharp given the marginal distributions of potential outcomes. 
In the presence of covariates, 
%In the current setting considered here, the bounds given by \citet{neyman1990application} and \citet{Aronow2014SBotViRE} are valid but does not make use of the information contained in $w_i$, $i=1,2,\dots,N$. 
%When $w_i=(1,x_{i1},\dots,x_{ip})^T$ is a vector of covariates, 
\citet{Ding2018DTEV} proposed the following regression based bound that may   improve \cite{Aronow2014SBotViRE}'s bound in certain situations:
\begin{align*}
	\phi^2_{\rm DL}\colonequals&\phi^2(\tau_{w})+\int_0^{1}(F_{e_1}^{-1}(u)-F_{e_0}^{-1}(u))^2du\leq\phi^2(\tau)
	\\&\leq\phi^2(\tau_{w})+\int_0^{1}(F_{e_1}^{-1}(u)-F_{e_0}^{-1}(1-u))^2du\equalscolon\phi^2_{\rm DH}.
\end{align*}
where $\tau_{w}=(w_1^\T(\gamma_1-\gamma_0),\dots,w_N^\T(\gamma_1-\gamma_0))^\T$, $F_{e_t}(s)=P(e_t\leq s)$,  $e_{t}=(y_{t1}-w_1^\T\gamma_t,\dots,y_{tN}-w_N^\T\gamma_t)^\T$,  and $\gamma_t$ is the least square regression coefficient of $y_{ti}$ on $w_i$.
The lower bound of \cite{Ding2018DTEV}  is not sharp, as we observe that it can even be smaller than  that of \cite{Aronow2014SBotViRE} and thus may lead to more conservative confidence intervals in spite of the available covariate information. Such situation does not occurs with our bound. It can be verified that the bounds $\phi_{\rm AL}^{2}$, $\phi_{\rm AH}^{2}$, $\phi_{\rm DL}^{2}$, $\phi_{\rm DH}^{2}$ are all functionals of $\pi_{k}$ and $F_{t\mid k}$ ($t=0,1$ and $k=1,\dots, K$). Thus we have
%We can see this from Fig \ref{fig: ex1}.
\begin{equation}\label{eq:bounds relationship}
	\phi^2_{\rm L}\geq\max\{\phi^2_{\rm AL},\phi^2_{\rm DL}\}, \ \phi^2_{\rm H}\leq\min\{\phi^2_{\rm AH},\phi^2_{\rm DH}\}
\end{equation}
according to Theorem \ref{thm:cov}.
When there is no covariate,  our bound reduces to $[\phi^2_{\rm AL},  \phi^2_{\rm AH}]$ by letting  $K=1$, $\xi_{1}=1$ and $w_{i}=1$ for $i=1,\dots,N$. 
%From the structure of $\phi^2_{\rm AL}$ and $\phi^2_{\rm DL}$, we can see that $\phi^2_{\rm AL}$ and $\phi^2_{\rm DL}$ are actually lower bounds of the set $\{\phi^2(\tau^*):\tau^*=y_1^*-y_0^*,\mathbf{U}^*=(y_1^*,y_0^*,w^*) \in \mathcal{U}\}$ and our bound $\phi_{\rm L}^2$ can be attained by some population in $\mathcal{U}$. 
The following example  illustrates the improvement of our bound as the association between covariates and potential outcomes varies. 
\begin{example}\label{ex: bound perform}
	Consider a population with $N=600$ units. Suppose the potential outcomes and the   covariate  are  binary, 
	with $P(w=1)=1/3$, $P(y_1=1)=2/3$, $P(y_0=1)=1/3$, $P(y_0=1 \mid w=1)=3/4$, and $P(y_0=1 \mid w=0)=1/8$. 
	Let $p=P(y_1=1\mid w=1)$ ($p\in\{1/200,\dots,1\}$), then $P(y_1=1\mid w=0)=1-p/2$. 
	Figure 1 presents the three bounds   under different values of $p$.
\end{example}

\begin{figure}[h]
	% The arguments in the next line are {height}{optional width}{used only by OUP for typesetting}[filename, in directory art]
	\centering
	\subfigure[Relationship between $\phi^2_{\rm AL}$, $\phi^2_{\rm DL}$ and $\phi^2_{\rm L}$.]{
		\begin{minipage}[t]{0.45\textwidth}
			\centering
			\includegraphics[scale=0.4]{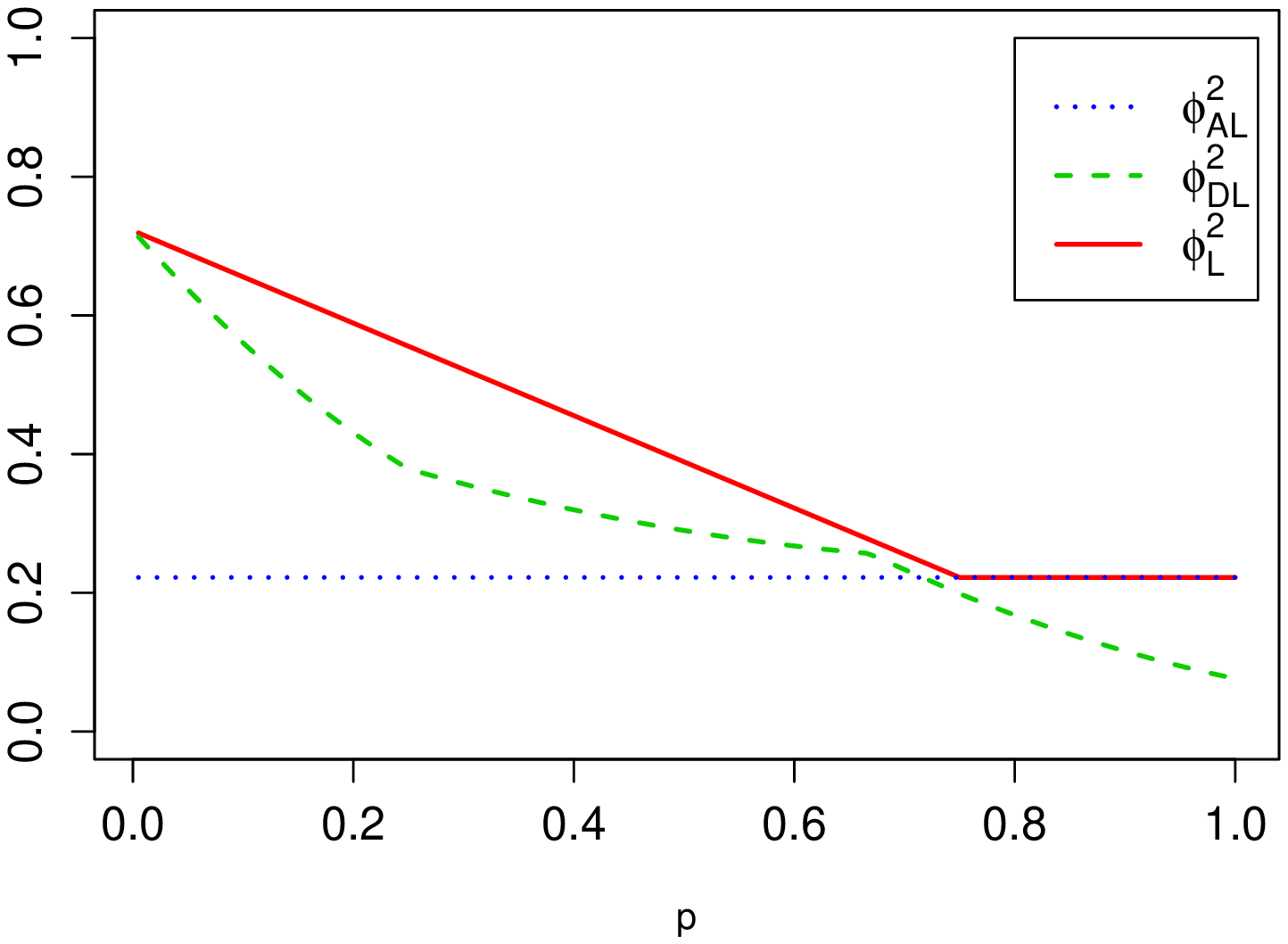}
		\end{minipage}
	}
	\subfigure[Relationship between $\phi^2_{\rm AH}$, $\phi^2_{\rm DH}$ and $\phi^2_{\rm H}$.]{  
		\begin{minipage}[t]{0.45\textwidth}
			\centering
			\includegraphics[scale=0.4]{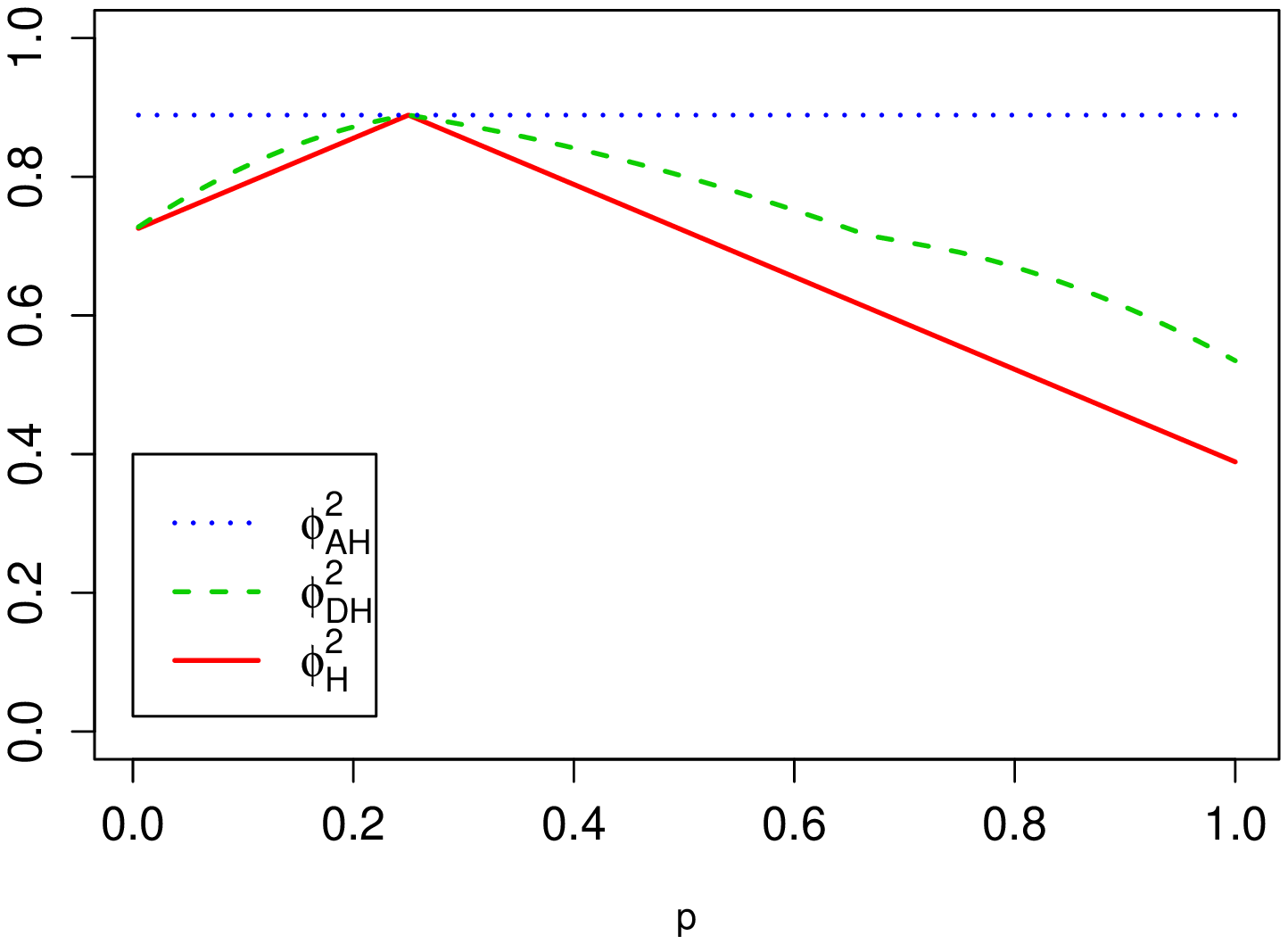}
		\end{minipage}
	}
	\caption{Comparison of three bounds under different values of $p$.}\label{fig: ex1}
\end{figure}

Figure 1 shows that $\phi^2_{\rm L}\geq\max\{\phi^2_{\rm AL},\phi^2_{\rm DL}\}$ and $\phi^2_{\rm H}\leq\min\{\phi^2_{\rm AH},\phi^2_{\rm DH}\}$ under all settings of $p$, 
and in many situations  the  inequalities are strict.  
For  $p\leq 143/200$, \cite{Ding2018DTEV}'s bound is tighter than \cite{Aronow2014SBotViRE}'s bound.
However, for $p>143/200$, $\phi^2_{\rm DL}<\phi^2_{\rm AL}$,  although covariate information is also used in the approach of \cite{Ding2018DTEV}.  
%By Cauchy-Schwartz inequality, we have $\max\{\phi^2(\tau) - \phi^2_{\rm DL}, \phi^2_{\rm DH} - \phi^2(\tau)\} \leq \phi^2_{\rm DH} - \phi^2_{\rm DL} \leq 4(\sum_{i=1}^Ne_{1i}^2/N)^{1/2}(\sum_{i=0}^Ne_{0i}^2/N)^{1/2}$. Thus \cite{Ding2018DTEV}'s lower bound may perform well when the residuals $e_{ti}$ are small, that is, a linear combination of the covariates fit the outcomes well. When the linear model does not fit the outcomes well, its performance is not guaranteed. However, the bound $[\phi^2_{\rm L}, \phi^2_{\rm H}]$ is always the best among the three bounds.

\subsection{Estimation of the sharp bound and the confidence interval}\label{subsec: estimation cov}
To estimate $\phi_{\rm L}^2$ and $\phi_{\rm H}^2$ and study asymptotic properties of the
proposed estimators for them, we adopt the following standard framework \citep{Li2017GFoFPCLTwAtCI} for theoretical development. Suppose that there is a sequence of finite populations
$\mathbf{U}_{N}$ of size $N$. For each $N$, $n_{1}$ units are randomly assigned to treatment group and $n_{0}$ units to control group. As the population size $N\to\infty$, the sizes of the treatment and the control groups satisfy $n_{1}/N\to\rho_1$, $n_{0}/N\to\rho_0$ with $\rho_1$, $\rho_0\in(0,1)$ and $\rho_1+\rho_0\leq 1$. Here we assume that the number of covariate values $K$ is known and is  allowed to grow at a certain rate with the population size $N$.
%\begin{enumerate}
%	\item $K$ is fixed while $N\to \infty$;
%	\item $K\to \infty$ as $N\to \infty$.
%\end{enumerate}
%Scenario 1 holds  if all  covariates  are categorical or ordinal variables. 
To accommodate   continuous covariates, we can stratify them and increase  the number of strata with the sample size. 
We estimate $\pi_k$ and $F_{t\mid k}(y)$ with  empirical probabilities $\hat{\pi}_k=1/n\sum_{T_i\in \{0,1\}}1\{w_i=\xi_k\}$ and $\hat{F}_{t\mid k}(y)=\sum_{T_{i} = t}1\{y_{ti}\leq y, w=\xi_k\}/\sum_{T_{i} = t}1\{w_{i} = \xi_{k}\}$, respectively. 
By plugging  in these estimators, we obtain the following estimators for $\phi_{\rm L}^{2}$ and $\phi_{\rm H}^{2}$:
\begin{equation}\label{eq: est cov}
	\begin{aligned}
		&\hat \phi^2_{\rm L} = \sum_{k=1}^{K}\hat{\pi}_k\int_0^{1}(\hat{F}_{1\mid k}^{-1}(u)-\hat{F}_{0\mid k}^{-1}(u))^2du-\hat{\theta}^2,\\
		&\hat \phi^2_{\rm H} =  \sum_{k=1}^{K}\hat{\pi}_k\int_0^{1}(\hat{F}_{1\mid k}^{-1}(1-u)-\hat{F}_{0\mid k}^{-1}(u))^2du-\hat{\theta}^2.
	\end{aligned}
\end{equation}
These estimators involve sample conditional quantile functions $\hat{F}_{t\mid k}^{-1}(u)$, whose statistical properties are complicated in the finite population framework. Many theoretical results for conditional quantile functions in the super population framework can not be applied to the scenario we consider. Moreover, the number of covariate values $K$ is allowed to diverge as the population size, which further complicates the problem. Thus it is not trivial to analyze the asymptotic properties of these estimators. However, we observe that the first term of these estimators is actually weighted sums of the Wasserstein distances between some distributions. By invoking a representation theorem of Wasserstein distances, we are able to prove the consistency results for the estimators with a careful analysis of the error terms. See the Appendix for more details.
We assume that the population $\mathbf{U}_{N}$ satisfies the following two conditions.
\begin{condition}\label{cond:fourth moment}
	There is some constant $C_M$ which does not depend on $N$ such that $1/N\sum_{i=1}^{N}y_{ti}^4\leq C_M$ for $t=0,1$.
\end{condition}
\begin{condition}\label{cond:positivity}
	There is some constant $C_{\pi}$ which does not depend on $N$ such that $\pi_{k}\geq C_{\pi}/K$ for $k=1,\dots,K$.
\end{condition}
Furthermore, we assume $K$ satisfies the following condition.
\begin{condition}\label{cond:order1}
	$K^2\log K/N \to 0$ as $N \to \infty$.
\end{condition}
Condition \ref{cond:fourth moment} requires the potential outcomes have uniformly bounded fourth moments. Condition \ref{cond:positivity} requires the proportion of units with each covariate value not to be too small. Condition \ref{cond:order1} imposes some upper bound on the number of values the covariate may take. If the covariate $w$ is some subgroup indicator, then Condition \ref{cond:order1} can be easily satisfied if the number of subgroups is not too large. Continuous components in $w$ can be stratified to meet Condition \ref{cond:order1}. If $w$ contains many components, then the number of covariate values may be too large even after stratifying the continuous components. In this case, one can partition units with similar covariate values into subgroups and use the subgroup label as the new covariate to apply the proposed method. Alternatively, in practice, one can first conduct dimension reduction or variable selection methods to obtain a low-dimensional covariate, and then apply the proposed method using the obtained covariate.
The following theorem establishes the consistency of the bound estimators.
\begin{theorem}\label{thm:consistency cov}
	Under Conditions \ref{cond:fourth moment}, \ref{cond:positivity} and \ref{cond:order1}, we have  
	\[
	(\hat{\phi}_{\rm L}^2, \hat{\phi}_{\rm H}^2) - (\phi_{\rm L}^2, \phi_{\rm H}^2) \stackrel{P}{\to} 0.
	\]
	as $N\to \infty$.
\end{theorem}

Proof of this theorem is relegated to the Appendix.
The lower bound $\phi^{2}_{\rm L}$ for $\phi^{2}(\tau)$ implies an upper bound for $\sigma^{2}$. 
Consistency result of $\hat{\phi}^{2}_{\rm L}$ is sufficient for constructing a conservative confidence interval.
A conservative $1-\alpha$ confidence interval for $\theta$ is given by
\begin{equation}\label{eq: ci ate}
	I_{N} = \left[\hat{\theta} - q_{\frac{\alpha}{2}}\hat{\sigma}N^{-\frac{1}{2}}, \hat{\theta} + q_{\frac{\alpha}{2}}\hat{\sigma}N^{-\frac{1}{2}}\right],
\end{equation}
where \[\hat{\sigma}^{2} = \left(\frac{N}{n_1}\hat{\phi}^2_{1}+\frac{N}{n_0}\hat{\phi}^2_{0}-\hat{\phi}^2_{L}\right)\]
and $q_{\alpha/2}$ is the upper $\alpha/2$ quantile of a standard normal distribution.

Next, we study the property of the confidence interval $I_{N}$ in \eqref{eq: ci ate}. As discussed in \citep{lehmann2006testing}, inference based on asymptotic results is not reassuring unless some uniform convergence results can be established. 
% The asymptotic results characterize the performance of statistical methods if the sample size is sufficiently large. Then a natural question is ``how large is sufficiently large in a specific problem?". If no uniform convergence result is available, ``how large is sufficient large" may depend on the data distribution which is unknown. So we can not ensure whether the sample size is sufficiently large. 
We next show that $I_{N}$ is uniformly asymptotically level $1-\alpha$  over a class of finite populations. For some constants $L_{1}, L_{2}, L_{3} > 0$, we introduce the following class of finite populations
\begin{equation}\label{eq: def pop class}
	\begin{aligned}
		\mathcal{P}_{N} & = \Bigg\{\mathbf{U}_{N}^{*}=(y_1^*,y_0^*,w^*):\ \text{$\mathbf{U}_{N}^{*}$ is of size $N$, and}\\
		&\qquad \text{(a)}\ \frac{1}{N} \sum_{i=1}^{N} y_{ti}^{*4}\leq L_{1}; \quad \text{(b)} \
		\phi^{2}(y_{t}^{*}) \geq L_{2}; \quad \text{(c)}\ P(w^*=\xi_k)\geq \frac{L_{3}}{K}\\
		&\qquad \text{for $t=0,1$ and $k=1,\dots,K$}\Bigg\} .
	\end{aligned}
\end{equation}
Under Conditions \ref{cond:fourth moment} and \ref{cond:positivity}, if the variances of the potential outcomes are bounded away from zero, then $\mathbf{U}_{N}$ belongs to $\mathcal{P}_{N}$ for some $L_{1}$, $L_{2}$ and $L_{3}$.
Constraint (a) in the definition of $\mathcal{P}_{N}$ requires the fourth moments of the potential outcomes to be uniformly bounded. Bounded fourth moments are required for theoretical development in many existing works \citep{Freedman2008ORAiEwST,freedman2008regression}. Constraint (b) requires the variance of the potential outcomes to be bounded away from zero. According to Cauchy-Schwartz inequality and some straightforward calculations, Constraint (b) implies that the variance of $\sqrt{N}(\hat{\theta} - \theta)$ is bounded away from zero at least for sufficiently large $N$ if $\mathbf{U}_{N} \in \mathcal{P}_{N}$. Constraint (c) requires units with each covariate value not to be too rare. Constraints (b) and (c) ensure that the the denominator of some quantities appears in the theoretical analysis are not too small, which is important to establish the desired properties. The set $\mathcal{P}_{N}$ contains a large class of finite populations. For an illustration, suppose $(y_{1i}, y_{0i}, w_{i})$, $i=1,\dots,n$, are independent and identically distributed observations of some random variables $(Y_{1}, Y_{0}, W)$. Then according to the strong law of large number, $\mathbf{U}_{N}$ belongs to $\mathcal{P}_{N}$ for sufficiently large $N$ with probability one as long as $(Y_{1}, Y_{0}, W)$ satisfies $E[Y_{t}^{4}] < L_{1}$, $\var[Y_{t}] > L_{2}$ and $\mathbb{P}(W = \xi_{k}) > L_{3}/K$ for $t = 0,1$ and $k = 1,\dots,K$.
%\begin{equation}\label{eq: dist class}
%	\begin{aligned}
	%		\Bigg\{
	%		&\mathbb{P}:\ \text{under} \ \mathbb{P},\ \E(Y_{t}^{4}) \leq \bar{L}_{1}, \ \var(Y_{t})\geq \bar{L}_{2}, \mathbb{P}(W = \xi_{k}) \geq \frac{\bar{L}_{3}}{K} \\
	%		&\text{for $t=0,1$ and $k=1,\dots,K$}
	%		\Bigg\}.
	%	\end{aligned}
%\end{equation}
%Here $\bar{L}_{1}$, $\bar{L}_{2}$, $\bar{L}_{3}$ is some constants satisfying $\bar{L}_{1} < L_{1}$, $\bar{L}_{2} > L_{2}$ and $\bar{L}_{3} > L_{3}$.
%%Without loss of generality, we assume $\mathcal{P}_{N}$ is non-empty at least for sufficiently large $N$. For $N \geq 2$, 
%The class \eqref{eq: dist class} is a large nonparametric class of distributions.

%For example, $L_{1} \geq 4 L_{2}^{2}(\rho_{1}^{-1} + \rho_{2}^{-1})^{-2}$ and $L_{3} \leq 1$ can ensure that $\mathcal{P}_{N}$ is non-empty. 
Next, we show that $I_{N}$ has uniformly asymptotically guaranteed coverage rate over $\mathcal{P}_{N}$ for any $L_{1}, L_{2}, L_{3} > 0$.
\begin{theorem}\label{thm: asy level cov}
	Under Condition \ref{cond:order1}, the confidence interval $I_{N}$ in \eqref{eq: ci ate} is uniformly asymptotically level $1 - \alpha$ over $\mathcal{P}_{N}$ for any $L_{1}, L_{2}, L_{3} > 0$, that is
	\begin{equation*}
		\mathop{\lim\inf}_{N\to \infty} \inf_{\mathbf{U}_{N}\in\mathcal{P}_{N}} \mathbb{P}(\theta \in I_{N}) \geq 1-\alpha.
	\end{equation*}
\end{theorem}

Proof of this theorem is in the Appendix. 

%%%%%%%%%%%%%%%%%%%%%%%%%%%%%%%%%%%%%%%%%%%%%%%%%%%%%%%%%%%%%%%%%%%%%%%%%%%%%%%%%%%%%%%%%%%%%%%%%%%%%%%%%%%%%%%%%%%%%%%%%%%%

\section{Sharp variance bound for the Wald estimator}\label{sec:noncom}
\subsection{Sharp bound for the unidentifiable term in the variance}
In the previous section, we discuss the variance bound in completely randomized experiments with perfect compliance where each unit takes the treatment assigned by the randomization procedure. However, noncompliance often arises in randomized experiments. In the presence of noncompliance, some units may take the treatment different from the assigned treatment following their own will or for other reasons. For each unit $i$ and $t=0,1$, we let $d_{ti}\in\{0,1\}$ denote the treatment that unit $i$ actually takes if assigned to treatment $t$. In this case, the units can be classified into four groups according to the value of $(d_{1i},d_{0i})$ \citep{angrist1996identification,frangakis2002principal},
\[
g_i=\left\{\begin{array}{lcl}
	\text{Alway Taker}\ (\rm a)& \qquad &\text{if} \ d_{1i}=1\ \text{and} \ d_{0i}=1,\\
	\text{Complier}\ (\rm c)& \qquad &\text{if} \ d_{1i}=1\ \text{and} \ d_{0i}=0,\\
	\text{Never Taker}\ (\rm n)& \qquad &\text{if} \ d_{1i}=0\ \text{and} \ d_{0i}=0,\\
	\text{Defier}\ (\rm d)& \qquad &\text{if} \ d_{1i}=0\ \text{and} \ d_{0i}=1.\\
\end{array}
\right.
\]
Let $g=(g_1,\dots,g_N)$. Then the characteristics of the population  can be viewed as a matrix $\mathbf{U}_{\rm c}=(y_1,y_0,w,g)$ where $y_1$, $y_0$ and $w$ are defined in Section \ref{sec:CRT}.
For $t=0,1$, $k=1,\dots,K$ and $h=\rm a,c,n, d$, let $F_{t\mid (k,h)}(y) = P(y_t \leq y\mid w = \xi_{k}, g =  h)$, $\pi_{k\mid h} = P(w = \xi_k\mid g= h)$ and $\pi_{h}=P(g=h)$.

In this section, we maintain the following standard assumptions for analyzing the randomized experiment with noncompliance.

\begin{assumption}\label{ass:iv}(i) Monotonicity:$d_{1i}\geq d_{0i}$; (ii) exclusion restriction: $y_{1i}=y_{0i}$ if $d_{1i}=d_{0i}$; (iii) strong instrument: $\pi_{\rm c}\geq C_{0}$ where $C_{0}$ is a positive constant.
\end{assumption}
Assumption \ref{ass:iv} (i) rules out the existence of defiers and is usually easy to assess in many situations. For example, it holds automatically if units in the control group do not have access to the treatment.
Assumption \ref{ass:iv} (ii) means that the treatment assignment affect the potential outcome only through affecting the treatment a unit actually receive. Assumption \ref{ass:iv} (iii) assures the existence of compliers. Assumption \ref{ass:iv} is commonly adopted to identify the causal effect in the presence of noncompliance. See \cite{angrist1996identification,abadie2003semiparametric} for more detailed discussions of Assumption \ref{ass:iv}. 
In the randomized experiment with noncompliance, the parameter of interest is the local average treatment effect (LATE) \citep{angrist1996identification,abadie2003semiparametric},
\[
\theta_{\rm c} = \frac{\sum_{i=1}^{N}\tau_{i}1\{g_i={\rm c}\}}{\sum_{i=1}^{N}1\{g_i={\rm c}\}},
\]
which is the average effect of treatment in the compliers.

Under the monotonicity and exclusion restriction, we have $1\{g_i ={\rm c}\}= d_{1i} - d_{0i}$ and $(d_{1i} - d_{0i})\tau_i = \tau_i=y_{1i}-y_{0i}$. Thus
\[\pi_{\rm c} = \frac{1}{N}\sum_{i=1}^N1\{g_i= {\rm c}\} = \frac{1}{N}\sum_{i=1}^{N}(d_{1i}-d_{0i}) = \mu(d_1) - \mu(d_0),\]
\[\frac{1}{N}\sum_{i=1}^{N}(d_{1i} - d_{0i})\tau_i = \frac{1}{N}\sum_{i=1}^{N}\tau_{i} = \mu(\tau),\]
and $\theta_{\rm c} = \pi_{\rm c}^{-1}\theta$.
Hence $\theta_{\rm c}$ can be estimated by the ``Wald estimator",
\[\hat{\theta}_{\rm c} = \hat{\pi}_{\rm c}^{-1}\hat{\theta}\]
where $\hat{\theta} = \sum_{T_i=1}y_{1i}/n_{1} - \sum_{T_i=0}y_{0i}/n_{0}$ and
$\hat{\pi}_{\rm c} = \sum_{T_i=1}d_{1i}/n_{1} - \sum_{T_i=0}d_{0i}/n_{0}$.
Let $z_{i} = (y_{1i}, y_{0i}, d_{1i}, d_{0i})^{\T}$, $\bar{z} = \sum_{i=1}^{N}z_i/N$ and 
\[V_{N} = \frac{1}{N}\sum_{i=1}^{N}(z_{i}-\bar{z})(z_{i}-\bar{z})^{\T}.\] 
The asymptotic normality is established under the following regularity condition.
\begin{condition}\label{cond:nonsingular cov}
	There is some constant $C_{\lambda}$  which does not dependent on $N$ such that the eigenvalues of $V_{N}$
	is not smaller than $C_{\lambda}$.
\end{condition}
Now we are ready to state the asymptotic normality result.
\begin{theorem}\label{thm:asymptotic dist compliance}
	Under Assumption \ref{ass:iv} and Conditions \ref{cond:fourth moment}, \ref{cond:nonsingular cov}, we have
	\[\sqrt{N}\sigma_{\rm c}^{-1}(\hat{\theta}_{\rm c} - \theta_{\rm c})\overset{d}{\to}N(0, 1)\]
	provided $\sigma^{2}_{\rm c}$ is bounded away from zero,
	where 
	\[\sigma^{2}_{\rm c} = \frac{1}{\pi_{\rm c}^{2}}\left(\frac{N}{n_1}\phi^2(\tilde{y}_{1})+\frac{N}{n_0}\phi^2(\tilde{y}_{0})-\phi^2(\tilde{\tau})\right),\]
	$\tilde{y}_{t} = (y_{t1} - \theta_{\rm c} d_{t1},\dots, y_{tN} - \theta_{\rm c} d_{tN})^\T$ for $t=0,1$ and $\tilde{\tau} = \tilde{y}_{1} - \tilde{y}_{0}$.
\end{theorem}

Proof of this theorem is in the Appendix.
Let $\hat{y}_{ti} = y_{ti} - \hat{\theta}_{\rm c}d_{ti}$, then under the conditions of Theorem \ref{thm:asymptotic dist compliance}, $\phi^2(\tilde{y}_{t})$ can be consistently estimated by 
\begin{equation}\label{eq: est phi_check_t}
	\check{\phi}^2_{t} =  \frac{1}{n_t - 1}\sum_{T_i=t}\Big(\hat{y}_{ti} - \frac{1}{n_t}\sum_{T_j=t}\hat{y}_{tj}\Big)^2,
\end{equation}
and $\pi_{\rm c}$ can be consistently estimated by $\hat{\pi}_{\rm c}$. However, analogous to $\phi^2(\tau)$, $\phi^2(\tilde{\tau})$ is unidentifiable. Here we construct a sharp bound for $\phi^2(\tilde{\tau})$ using covariate information. It should be pointed out that the sharp 
bound has not been obtained even in the absence of  covariates. Define $\tilde{F}_{t\mid k}(y) = P(\tilde{y}_{t}\leq y\mid w=\xi_{k},g={\rm c})$ be a set of lower bounds for $\phi^{2}(\tilde{\tau})$. Let $\tilde{\mathcal{B}}_{\rm L} = \{\tilde{b}_{\rm L}: \text{$\tilde{b}_{\rm L}$ is a functional of $\pi_{\rm c}$, $\pi_{k\mid \rm c}$ and $\tilde{F}_{t\mid k}$ for $t=0,1$ and $k=1,\dots, K$}; \tilde{b}_{\rm L} \leq \phi^{2}(\tilde{\tau})\}$. Define the set of upper bounds $\tilde{\mathcal{B}}_{\rm H}$ similarly. Then we can establish the following sharp bound for $\phi^2(\tilde{\tau})$.
\begin{theorem}\label{thm:compliance}
	A bound for  $\phi^2(\tilde{\tau})$ is $[\tilde{\phi}^2_{\rm L}, \tilde{\phi}^2_{\rm H}]$ where
	\begin{align*}
		&\tilde{\phi}^2_{\rm L}=\sum_{k=1}^{K}\pi_{\rm c}\pi_{k\mid \rm c}\int_0^{1}(\tilde{F}_{1\mid k}^{-1}(u)-\tilde{F}_{0\mid k}^{-1}(u))^2du,\\
		&\tilde{\phi}^2_{\rm H}=\sum_{k=1}^{K}\pi_{\rm c}\pi_{k\mid \rm c}\int_0^{1}(\tilde{F}_{1\mid k}^{-1}(u)-\tilde{F}_{0\mid k}^{-1}(1-u))^2du.
	\end{align*}
	Moreover, the bound is sharp in the sense that $\tilde{\phi}_{\rm L}^{2}$ is the largest lower bound in $\tilde{\mathcal{B}}_{\rm L}$ and  $\tilde{\phi}_{\rm H}^{2}$ is the smallest upper bounds in $\tilde{\mathcal{B}}_{\rm H}$.
\end{theorem}

See the Appendix for the proof of this theorem.
If there is no covariate, we let  $K=1$, $\xi_{1}=1$ and define $w_{i}=1$ for $i=1,\dots,N$, then we can obtain a bound without covariates, which has not been considered in the literature.

\subsection{Estimation of the sharp bound and the confidence interval}
To estimate the bounds, we need to estimate $\pi_{\rm c}$, $\pi_{k\mid \rm c}$, and
$\tilde{F}_{t\mid k}(y)$. Here $\pi_{\rm c}$ can be estimated by $\hat{\pi}_{\rm c}$.
We let
%According to the monotonicity and exclusion restriction, we have $1\{g_i = {\rm c}\} = d_{1i} - d_{0i}$, $(1 - d_{1i})1\{\tilde{y}_{0i}\leq y\} = (1 - d_{1i})1\{\tilde{y}_{1i}\leq y\}$ and $d_{0i}1\{\tilde{y}_{1i}\leq y\} = d_{0i}1\{\tilde{y}_{0i}\leq y\}$. Thus,
%\begin{align*}
%\pi_{k\mid{\rm c}} 
%= &\pi_{\rm c}^{-1} \left(\frac{1}{N}\sum_{i=1}^{N}d_{1i}1\{w_{i}=\xi_{k}\} - \frac{1}{N}\sum_{i=1}^{N}d_{0i}1\{w_{i}=\xi_{k}\}\right)\\
%= &\pi_{\rm c}^{-1} \left(\frac{1}{N}\sum_{i=1}^{N}(1 - d_{0i})1\{w_{i}=\xi_{k}\} - \frac{1}{N}\sum_{i=1}^{N}(1 - d_{1i})1\{w_{i}=\xi_{k}\}\right),
%\end{align*}
%\begin{equation*}
%\tilde{F}_{1\mid k}(y) 
%= \pi_{\rm c}^{-1}\pi_{k\mid{\rm c}}^{-1}\left(\frac{1}{N}\sum_{i=1}^{N}d_{1i}1\{\tilde{y}_{1i}\leq y\}1\{w_{i}=\xi_{k}\} - \frac{1}{N}\sum_{i=1}^{N}d_{0i}1\{\tilde{y}_{0i}\leq y\}1\{w_{i}=\xi_{k}\}\right),
%\end{equation*}
%and
%\begin{equation*}
%\tilde{F}_{0\mid k}(y) = \pi_{\rm c}^{-1}\pi_{k\mid{\rm c}}^{-1}\left(\frac{1}{N}\sum_{i=1}^{N}(1-d_{0i})1\{\tilde{y}_{0i}\leq y\}1\{w_{i}=\xi_{k}\} - \frac{1}{N}\sum_{i=1}^{N}(1 - d_{1i})1\{\tilde{y}_{1i}\leq y\}1\{w_{i}=\xi_{k}\}\right).
%\end{equation*}
%
%Based on these relationships, we define
\begin{align*}
	&\hat{\lambda}_{1k} = \frac{1}{n_1}\sum_{T_i=1}d_{1i}1\{w_{i}=\xi_{k}\} - \frac{1}{n_0}\sum_{T_i=0}d_{0i}1\{w_{i}=\xi_{k}\} \\
	&\hat{\lambda}_{0k} = \frac{1}{n_0}\sum_{T_i=0}(1 - d_{0i})1\{w_{i}=\xi_{k}\} - \frac{1}{n_1}\sum_{T_i=1}(1 - d_{1i})1\{w_{i}=\xi_{k}\}.
\end{align*}
Under Assumption \ref{ass:iv}, we estimate $\pi_{k\mid \rm c}$  and $\tilde{F}_{t\mid k}(y)$ with 
\[
\hat{\pi}_{k \mid {\rm c}} = \hat{\pi}_{\rm c}^{-1}\hat{\lambda}_{1\mid k},
\]
\[
\begin{aligned}
	\check{F}_{1\mid k}(y)& =\hat{\lambda}_{1\mid k}^{-1}\Bigg(\frac{1}{n_1}\sum_{T_i=1}d_{1i}1\{\hat{y}_{1i}\leq y\}1\{w_{i}=\xi_{k}\}\\
	&\quad - \frac{1}{n_0}\sum_{T_i=0}d_{0i}1\{\hat{y}_{0i}\leq y\}1\{w_{i}=\xi_{k}\}\Bigg)
\end{aligned}
\]
and
\[\begin{aligned}
	\check{F}_{0\mid k}(y)& = \hat{\lambda}_{0\mid k}^{-1}\Bigg(\frac{1}{n_0}\sum_{T_i=0}(1-d_{0i})1\{\hat{y}_{0i}\leq y\}1\{w_{i}=\xi_{k}\}\\ 
	&\quad - \frac{1}{n_1}\sum_{T_i=1}(1 - d_{1i})1\{\hat{y}_{1i}\leq y\}1\{w_{i}=\xi_{k}\}\Bigg)
\end{aligned}
\]
where $\hat{y}_{ti} = y_{ti} - \hat{\theta}_{\rm c}d_{ti}$ for $t = 0,1$, and $i = 1,\dots,N$.

We then obtain estimators for $\tilde{\phi}_{\rm L}^{2}$ and $\phi_{\rm H}^{2}$ by plugging these estimators into the expressions of Theorem \ref{thm:compliance},
\begin{equation}\label{eq: est compliance}
	\begin{aligned}
		&\check \phi^2_{\rm L} = \sum_{k=1}^{K}\hat{\pi}_{\rm c}\hat{\pi}_{k\mid \rm c}\int_0^{1}(\check{F}_{1\mid k}^{-1}(u)-\check{F}_{0\mid k}^{-1}(u))^2du,\\
		&\check \phi^2_{\rm H} =  \sum_{k=1}^{K}\hat{\pi}_{\rm c}\hat{\pi}_{k\mid \rm c}\int_0^{1}(\check{F}_{1\mid k}^{-1}(u)-\check{F}_{0\mid k}^{-1}(1-u))^2du.
	\end{aligned}
\end{equation}
The estimator $\check{F}_{t\mid k}(y)$  may not be monotonic with respect to $y$, which brings about some difficulties for theoretical analysis. However, $\check{F}_{t\mid k}^{-1}(u)$ is still well defined, and in the next theorem we show that the non-monotonicity of $\check{F}_{t\mid k}(y)$ does not diminish the consistency of $\check \phi_{\rm L}^2$ and $\check \phi_{\rm H}^2$.

In the following asymptotic analysis, we denote $\mathbf{U}_{\rm c}$ by $\mathbf{U}_{{\rm c}, N}$. We assume that $\mathbf{U}_{{\rm c}, N}$ satisfies Assumption \ref{ass:iv}, Condition \ref{cond:fourth moment} and the following two conditions.
The following conditions are modified versions of Conditions \ref{cond:positivity} and \ref{cond:order1} in the presence of noncompliance.
\begin{condition}\label{cond:positivity-compliance}
	There is some constant $C_{\pi,{\rm c}}$ which does not depend on $N$ such that $\pi_{k\mid {\rm c}}\pi_{\rm c}\geq C_{\pi,{\rm c}}/K$ for $k=1,\dots,K$.
\end{condition}
\begin{condition}\label{cond:order2}
	$K^2\log K \max\{C_N^4, 1\}/N \to 0$ as $N \to \infty$,
	where $C_N = \max_{t,i}\vert y_{ti}\vert$.
\end{condition}
Then we are ready to establish the consistency of the estimators proposed in \eqref{eq: est compliance}.
\begin{theorem}\label{thm:consistency compliance}
	Under Assumption \ref{ass:iv} and Conditions \ref{cond:fourth moment}, \ref{cond:positivity-compliance}, \ref{cond:order2}, we have
	\[(\check \phi_{\rm L}^2, \check \phi_{\rm H}^2) - (\tilde\phi_{\rm L}^2, \tilde\phi_{\rm H}^2) \stackrel{P}{\to} 0.\] 
\end{theorem}

Proof of this theorem is relegated to the Appendix.
The lower bound $\tilde{\phi}^{2}_{\rm L}$ for $\phi^{2}(\tilde{\tau})$ implies an upper bound for $\sigma^{2}_{\rm c}$. 
By Theorems \ref{thm:asymptotic dist compliance} and \ref{thm:consistency compliance} we can construct a conservative $1-\alpha$ confidence interval for $\theta_{\rm c}$,
\begin{equation}\label{eq: ci late}
	I_{{\rm c}, N} = \left[\hat{\theta}_{\rm c} - q_{\frac{\alpha}{2}}\hat{\sigma}_{\rm c}N^{-\frac{1}{2}}, \hat{\theta}_{\rm c} + q_{\frac{\alpha}{2}}\hat{\sigma}_{\rm c}N^{-\frac{1}{2}}\right],
\end{equation}
where \[\hat{\sigma}_{\rm c}^{2} = \frac{1}{\hat{\pi}_{\rm c}^{2}}\left(\frac{N}{n_1}\check{\phi}^2_{1}+\frac{N}{n_0}\check{\phi}^2_{0}-\check{\phi}^2_{L}\right)\]
and $q_{\alpha/2}$ is the upper $\alpha/2$ quantile of a standard normal distribution. 
We then show that $I_{{\rm c},N}$ is uniformly asymptotically level $1-\alpha$  over a class of finite populations. In the following asymptotic analysis, we denote $\mathbf{U}_{\rm c}$ by $\mathbf{U}_{{\rm c}, N}$.
For any finite population $\mathbf{U}_{{\rm c}, N}^{*}=(y_1^*,y_0^*,w^*, g^{*})$, we define $\tilde{y}_{t}^{*}$ similarly as $\tilde{y}_{t}$ for $t=0,1$. Let $\Lambda_{N}^{*}$ be the smallest eigenvalue of 
\[
\frac{1}{N} \sum_{i=1}^{N}(z_{i}^{*} - \bar{z}^{*})(z_{i}^{*} - \bar{z}^{*})^{\T},
\] 
where $z_{i}^{*} = (y_{1i}^{*}, y_{0i}^{*}, d_{1i}^{*}, d_{0i}^{*})^{\T}$ and $\bar{z}^{*} = \sum_{i=1}^{N}z_{i}^{*}/N$.  For some constants $L_{0}, L_{1}, L_{2}, L_{3} > 0$, we introduce the following class of finite populations
\begin{equation*}
	\begin{aligned}
		\mathcal{P}_{{\rm c},N} & = \Bigg\{\mathbf{U}_{{\rm c}, N}^{*}=(y_1^*,y_0^*,w^*, g^{*}):\ \text{$\mathbf{U}_{{\rm c}, N}^{*}$ is of size $N$ and satisfies}\\
		&\qquad \text{(a)} \ \text{Assumption \ref{ass:iv};}\ \text{(b)}\ \Lambda_{N}^{*} \geq L_{0};\ 
		\text{(c)}\ \frac{1}{N} \sum_{i=1}^{N} y_{ti}^{*4}\leq L_{1};\\\ 
		&\qquad \text{(d)} \
		\phi^2(\tilde{y}_{t}^{*}) \geq L_{2}; \quad \text{(e)}\ P(w^*=\xi_k\mid g^{*} = {\rm c})P(g^{*} = {\rm c})\geq \frac{L_{3}}{K}\\
		&\qquad \text{for $t=0,1$ and $k=1,\dots,K$}\Bigg\}.
	\end{aligned}
\end{equation*}
Constraint (a) is required for the identification of LATE. Constraint (b) is a regularity condition to ensure the asymptotic normality of $\hat{\theta}_{\rm c}$. Constraints (c), (d) and (e) are similar to those constraints in the definition of $\mathcal{P}_{N}$ in Section \ref{subsec: estimation cov}. $\mathcal{P}_{{\rm c}, N}$ can be a large class of finite populations if $L_{0}$, $L_{2}$, $L_{3}$ are small and $L_{1}$ is large. The class of finite populations $\mathcal{P}_{{\rm c}, N}$ can be related to some class of generic distributions in the same way as discussed before Theorem \ref{thm: asy level cov}. The details are omitted here. Similar arguments as in the proof of Theorem \ref{thm: asy level cov} can show the following result.
\begin{theorem}\label{thm: asy level compliance}
	Under Condition \ref{cond:order2}, the confidence interval $I_{{\rm c}, N}$ in \eqref{eq: ci late} is uniformly asymptotically level $1 - \alpha$ over $\mathcal{P}_{{\rm c}, N}$, that is
	\begin{equation*}
		\mathop{\lim\inf}_{N\to \infty} \inf_{\mathbf{U}_{{\rm c}, N}\in\mathcal{P}_{{\rm c}, N}} \mathbb{P}(\theta \in I_{{\rm c}, N}) \geq 1-\alpha.
	\end{equation*}
\end{theorem}

\section{Simulations}\label{sec:simulation}
\subsection{Completely randomized experiments with perfect compliance}\label{subsec: sim cov}
In this subsection, we evaluate the performance of the  bounds and the estimators $\hat{\phi}_{\rm L}$, $\hat{\phi}_{\rm H}$ proposed in Section \ref{sec:CRT} via some simulations.  We first generate a finite population of size $N$ by drawing i.i.d. samples from the following data generation process:
\begin{enumerate}[(i)]
	\item $W$ takes value in $\{1,2,3,4\}$ with equal probability;
	\item for $w=1,2,3,4$, $Y_{1}\mid W=w\sim N(\mu_{w},\phi^2_{w})$, $V\mid W=w\sim N(0,6 - \phi^2_w)$, $Y_1\Perp V\mid W$ and $Y_0=0.3Y_1+V$ where $(\mu_1,\mu_2,\mu_3,\mu_4)=(3,0,-2,4)$ and $(\phi^2_1,\phi^2_2,\phi^2_3,\phi^2_4)=(2,1.5,5,4)$. 
\end{enumerate}
The generated values are viewed as the fixed finite population. We take $N = 400, 800$ and $2000$, respectively, to demonstrate the performance of the proposed method under different population sizes.
The following table exhibits the $\phi^{2}(\tau)$ and the true value of different bounds under different population sizes.

\begin{table}[h]
	\def~{\hphantom{0}}
	\centering
	\caption{Bounds $\phi^2_{\rm AL}$, $\phi^2_{\rm AH}$, $\phi^2_{\rm DL}$, $\phi^2_{\rm DH}$, $\phi^2_{\rm L}$, $\phi^2_{\rm H}$ and the true value of $\phi^{2}(\tau)$ under different population sizes}
	\vspace{10pt}
	\begin{tabular}{*{8}{c}}
		\toprule
		& $\phi^2_{\rm AL}$ & $\phi^2_{\rm AH}$ & $\phi^2_{\rm DL}$ &$\phi^2_{\rm DH}$ & $\phi^2_{\rm L}$ & $\phi^2_{\rm H}$ & $\phi^{2}(\tau)$ \\
		\midrule
		$N = 400$ & 0.98 & 58.14 & 1.02 & 58.04 & 9.04 & 42.92 & 16.72\\
		$N = 800$ & 0.78 & 54.79 & 0.74 & 54.76 & 8.83 & 40.53 & 17.70\\ 
		$N = 2000$ & 0.88 & 58.23 & 0.86 & 58.22 & 9.14 & 42.25 & 18.56\\
		\bottomrule
	\end{tabular}\label{table: true cov}
\end{table}

It can be seen that the intervals $[\phi_{\rm AL}^{2}, \phi_{\rm AH}^{2}]$, $[\phi_{\rm DL}^{2}, \phi_{\rm DH}^{2}]$, $[\phi_{\rm L}^{2}, \phi_{\rm H}^{2}]$ all contain $\phi^{2}(\tau)$ and hence all the bounds are valid. Under all population sizes, the lower bound $\phi_{\rm L}^{2}$ is much larger than the other two lower bounds, and  the upper bound $\phi_{\rm L}^{2}$ is much smaller than the other two upper bounds. 

The above results show the effectiveness of our bounds on the population level. Next, we consider the performance of the proposed bound estimators in completely randomized experiments.
In the simulation, half of the units are randomly assigned to the treatment group while the other half are assigned to the control group. Then we estimate the proposed bounds with the estimators defined in \eqref{eq: est cov}. The randomized assignment is repeated for $1000$ times. 
The root mean square error (RMSE) of the bound estimators under different $N$ is summarized in Table \ref{table: est cov}. The RMSE of the bound estimators decreases as the sample size increases, which confirms the consistency result in Theorem \ref{thm:consistency cov}.

\begin{table}[h]
	\def~{\hphantom{0}}
	\centering
	\caption{RMSE of the estimators for $\phi^2_{\rm L}$ and $\phi^2_{\rm H}$ under different population sizes ($n_{1} = n_{0} = N/2$)}
	\vspace{10pt}
	\begin{tabular}{*{4}{c}}
		\toprule
		& $N = 400$& $N = 800$ & $N = 2000$\\
		\midrule
		$\phi^2_{\rm L}$ & 1.71 & 0.87 & 0.66 \\
		$\phi^2_{\rm H}$ & 2.31 & 1.36 & 1.00\\
		\bottomrule
	\end{tabular}\label{table: est cov}
\end{table}

Next, we evaluate the performance of the bound estimators in constructing confidence intervals (CIs).
Different CIs can be constructed based on the asymptotic normality in \eqref{eq: asy-crt} and different lower bounds for $\phi^{2}(\tau)$. To obtain the CIs, we use $\hat{\phi}_{t}^{2}$ defined in \eqref{eq: est phi_hat_t} to estimate $\phi^2(y_t)$ for $t = 0,1$, and replace $\phi^{2}(\tau)$ in $\sigma^{2}$ by the estimators of different lower bounds. Bounds of \citet{Aronow2014SBotViRE} and \citet{Ding2018DTEV} are estimated by plug-in estimators as suggested in these works, and the proposed lower bound is estimated by the estimators defined in \eqref{eq: est cov}.
The average width (AW) and coverage rate (CR) of $95\%$ CIs based on the naive lower bound zero \citep{neyman1990application}, the estimator of $\phi^2_{\rm AL}$ \citep{Aronow2014SBotViRE}, the estimator of $\phi^2_{\rm DL}$ \citep{Ding2018DTEV} and the estimator of $\phi^2_{\rm L}$ are listed as follows:

\begin{table}[h]
	\def~{\hphantom{0}}
	\centering
	\caption{Average widths (AWs) and coverage rates (CRs) of 95\% CIs based on the naive bound, $\phi^2_{\rm AL}$, $\phi^2_{\rm DL}$ and $\phi^2_{\rm L}$ under different population sizes ($n_{1} = n_{0} = N/2$)}
	\vspace{10pt}
		\begin{tabular}{*{9}{c}}
			\toprule
			\multirow{2}{*}{Method}&\multicolumn{2}{c}{naive}&\multicolumn{2}{c}{$\phi^2_{\rm AL}$}&\multicolumn{2}{c}{$\phi^2_{\rm DL}$}
			&\multicolumn{2}{c}{$\phi^2_{\rm L}$}\\
			\rule{0pt}{19pt}
			&{\footnotesize AW} & {\footnotesize CR} &{\footnotesize AW} & {\footnotesize CR} &{\footnotesize AW} & {\footnotesize CR} &{\footnotesize AW} &{\footnotesize CR}\\
			\midrule		
			$N = 400$ & 1.511 & 98.0\% & 1.495 & 97.8\% & 1.493 & 97.8\% & 1.383 & 96.6\%\\
			$N = 800$ & 1.033 & 97.8\% & 1.025 & 97.7\% & 1.025 & 97.7\% & 0.945 & 96.7\%\\
			$N = 2000$ & 0.674 & 98.4\% & 0.668 & 98.3\% & 0.668 & 98.3\% & 0.619 & 97.3\%\\
			\bottomrule
		\end{tabular}
	\label{table: CI cov}
\end{table}

The CIs based on the estimator of $\phi_{\rm L}^{2}$ are the shortest, and the corresponding coverage rate is the closest to $95\%$ among the four CIs under all population sizes. See the Appendix for more simulation results on the performance of the proposed CI.

\subsection{Completely randomized experiments with noncompliance}\label{subsec: sim compliance}
Here we show the simulation performance of the  bounds and the estimators $\check{\phi}_{\rm L}$, $\check{\phi}_{\rm H}$ proposed in Section \ref{sec:noncom}. First, we
generate finite populations of size $N = 400, 800$ and $2000$ by drawing i.i.d. samples from the following data generation process:
\begin{enumerate}[(i)]
	\item generate the compliance type $G\in\{{\rm a, c, n}\}$ with the probability that $G = {\rm a}$, $\rm c$ and $\rm n$ being $0.2$, $0.7$ and $0.1$, respectively;
	\item for $h = \rm a,c$ and $\rm n$, the conditional distribution $W\mid G = h$ has probability mass $p_{1h}$, $p_{2h}$, $p_{3h}$ and $p_{4h}$ at $1,2,3,4$, respectively, where $(p_{1\rm a}, p_{2\rm a}, p_{3 \rm a}, p_{4\rm a}) = (0.15,0.2,0.3,0.35)$, $(p_{1\rm c}, p_{2\rm c}, p_{3 \rm c}, p_{4\rm c}) =$\\ $(0.25,0.25,0.25,0.25)$ and $(p_{1\rm n}, p_{2\rm n}, p_{3 \rm n}, p_{4\rm n}) = (0.35,0.3,0.2,0.15)$;
	\item for $w=1,2,3,4$, $Y_{1}\mid W=w\sim N(\mu_{w},\phi^2_{w})$ where $(\mu_1,\mu_2,\mu_3,\mu_4)=(3,0,-2,4)$ and $(\phi^2_1,\phi^2_2,\phi^2_3,\phi^2_4)=(2,1.5,5,4)$;
	\item
	$Y_{0}\mid G = {\rm c}, W=w \sim N(0.3w, 6 - \phi^2_w)$ and $Y_{0} = Y_{1}$
	if $G = \rm a$ or $\rm n$. 
\end{enumerate}
As stated behind Theorem \ref{thm:compliance}, the theorem can also provide a bound without using the covariate. So to illustrate the usefulness of the covariate, here we compare the bounds constructed using and without using the covariate. In the following, the lower bound using and without using the covariate is denoted by ``LC" and ``LNC", and the upper bound using and without using the covariate is denoted by ``HC" and ``HNC". As in Section \ref{subsec: sim cov}, half of the units are randomly assigned to the treatment group while the other half are  assigned to the control group. Then we estimate the bounds for each of these random assignments with the estimators proposed in \eqref{eq: est compliance}. The randomized assignment is repeated for $1000$ times. In the simulation $\phi^{2}(\tilde{\tau})$ equals to $9.97, 10.41$ and $9.56$ when $N = 400, 600$ and $2000$, respectively.
The following table exhibits the true values of different bounds and the RMSE of their estimators under different population sizes. 

\begin{table}[h]
	\def~{\hphantom{0}}
	\centering
	\caption{Bounds using and without using the covariate and the RMSE of their estimators under different population sizes. LNC: lower bound without covariate; HNC: upper bound without covariate; LC: lower bound with covariate; HC: upper bound with covariate ($n_{1} = n_{0} = N/2$)}
	\vspace{10pt}
		\begin{tabular}{*{9}{c}}
			\toprule
			& \multicolumn{2}{c}{LNC} & \multicolumn{2}{c}{HNC} &
			\multicolumn{2}{c}{LC} & \multicolumn{2}{c}{HC}\\
			& {\footnotesize Value} & {\footnotesize RMSE} & {\footnotesize Value} & {\footnotesize RMSE} & {\footnotesize Value} & {\footnotesize RMSE} & {\footnotesize Value} & {\footnotesize RMSE}\\
			\midrule
			
			$N = 400$ & 0.78 & 0.75 & 40.05 & 5.41 & 9.82 & 4.62 & 29.42 & 5.94\\
			$N = 800$ & 0.73 & 0.45 & 39.61 & 3.49 & 10.42 & 3.92 & 28.79 & 5.35\\
			$N = 2000$ & 0.72 & 0.27 & 37.19 & 2.10 & 9.52 & 2.45 & 28.45 & 4.55\\
			\bottomrule
		\end{tabular}\label{table: bound compliance}
\end{table}

Table \ref{table: bound compliance} shows that LC is much larger than LNC and HC is much smaller than HNC, which shows that the covariate is quite useful in sharpening the bounds. The RMSE of different bounds generally decreases as the sample size increases and this validates the consistency result in Theorem \ref{thm:consistency compliance}.

Next, we construct CIs based on the asymptotic normality in Theorem \ref{thm:asymptotic dist compliance} and different lower bounds for $\phi^{2}(\tilde{\tau})$.  To obtain the CIs, we use $\check{\phi}_{t}^{2}$ defined in \eqref{eq: est phi_check_t} to estimate $\phi^2(\tilde{y}_{t})$ for $t = 0,1$ and replace $\phi^{2}(\tilde{\tau})$ in $\sigma_{\rm c}^{2}$ by the estimators of different lower bounds.
The following table summarizes the average width (AW) and coverage rate (CR) of $95\%$ CIs based on the naive lower bound zero, and the estimator of lower bounds in Theorem \ref{thm:compliance} using and without using the covariate.

\begin{table}[h]
	\def~{\hphantom{0}}
	\centering
	\caption{Average widths (AWs) and coverage rates (CRs) of 95\% CIs based on the naive bound, HNL and HL under different population sizes.  LNC: lower bound without covariate; LC: lower bound with covariate ($n_{1} = n_{0} = N/2$)}
	\vspace{10pt}
		\begin{tabular}{*{7}{c}}
			\toprule
			\multirow{2}{*}{Method}&\multicolumn{2}{c}{naive}&\multicolumn{2}{c}{LNC}&\multicolumn{2}{c}{LC}\\
			&{\footnotesize AW} & {\footnotesize CR} &{\footnotesize AW} & {\footnotesize CR} &{\footnotesize AW} & {\footnotesize CR}\\
			\midrule
			$N = 400$ & 2.188 & 98.5\% & 2.168 & 98.4\% & 1.990 & 97.3\%\\
			$N = 800$ & 1.553 & 97.9\% & 1.542 & 97.8\% & 1.399 & 96.4\%\\
			$N = 2000$ & 0.980 & 98.2\% & 0.973 & 98.2\% & 0.893 & 96.6\%\\
			\bottomrule
		\end{tabular}
	\label{table: CI compliance}
\end{table}

The results show that the CIs based on the estimator of LC is the shortest, and the corresponding coverage rate is closest to $95\%$ among the three CIs under all population sizes. This demonstrates the usefulness of covariate in constructing CIs. 

\section{Real data applications}\label{sec:real data}
\subsection{Application to ACTG protocol 175}
In this section, we apply our approach proposed in Section \ref{sec:CRT} to a  dataset from the randomized trial ACTG protocol 175 \citep{hammer1996}  for illustration. Data used in this subsection are available from the R package ``speff2trial" (https://cran.r-project.org/web/ packages/speff2trial/index.\\html).
The ACTG 175 study evaluated four therapies in human immunodeficiency virus infected subjects whose CD4 cell counts (a measure of immunologic status) are from 200 to 500 $\rm mm^{-3}$. Here we regard the 2139 enrolled subjects as a finite population and consider two treatment arms: the standard zidovudine monotherapy (denoted by ``arm 0") and the combination therapy with zidovudine and didanosine (denoted by ``arm 1"). The parameter of interest is the average treatment effect of the combination therapy on the CD8 cell count measured at $20\pm 5$ weeks post baseline compared to the monotherapy. In the randomized trial, 532 subjects are randomly assigned to arm 0 and 522 subjects are randomly assigned to arm 1. The available covariate is the age of each subject. In order to meet Condition \ref{cond:order1},  we divide people into $\lfloor N^{1/4}\rfloor = 6$ groups according to age, with the first group less than 20 years old, the second group between 21 and 30 years old, the third group between 31 and 40 years old, the fourth group between 41 and 50 years old, the fifth group between 51 and 60 years old and the last group older than 60 years old. We then use the age group, gender and antiretroviral history as covariates and apply our method proposed in Section \ref{sec:CRT}. The proposed bounds $\phi^2_{\rm L}$ and $\phi^2_{\rm H}$ are estimated by the estimators defined in \eqref{eq: est cov}.
We also estimate the bounds proposed in the literature\citep{Aronow2013ACoUEotATEiRE,Ding2018DTEV}. Bounds $\phi^2_{\rm AL}$, $\phi^2_{\rm AH}$, $\phi^2_{\rm DL}$, $\phi^2_{\rm DH}$ are estimated by plug-in estimators as suggested in \cite{Aronow2014SBotViRE} and \cite{Ding2018DTEV}. The estimates of the lower bounds $\phi^2_{\rm AL}$, $\phi^2_{\rm DL}$ and $\phi^2_{\rm L}$ are $0.12$, $0.27$ and $4.75$, respectively; the estimates of the upper bounds $\phi^2_{\rm AH}$, $\phi^2_{\rm DH}$  and $\phi^2_{\rm H}$ are $70.79$, $69.87$ and $65.67$, respectively (values are divided by $10000$).
The estimate of the proposed lower bound is the largest among the three lower bounds and the estimate of the proposed upper bound is the smallest among the three upper bounds. The $95\%$ CI constructed using zero as a lower bound for $\phi^2(\tau)$ is $[-13.27, 92.79]$. Using the estimate of \citet{Aronow2014SBotViRE}'s lower bound lead to the CI $[-13.25, 92.77]$, and \citet{Ding2018DTEV}'s lower bound lead to the CI $[-13.23, 92.75]$. The CI $[-12.46, 91.98]$ is obtained by using $\hat{\phi}^{2}_{\rm L}$ given in \eqref{eq: est cov}. The widths of the four CIs are $106.07$, $106.02$, $105.98$ and $104.44$, respectively.
Comparing the CI width of the naive method using zero as the lower bound for $\phi^{2}(\tau)$  with the CI widths using the estimates of \cite{Aronow2013ACoUEotATEiRE}'s, \cite{Ding2018DTEV}'s and the proposed lower bound for $\phi^2(\tau)$, the CI width reductions are $0.04$, $0.09$ and $1.62$, respectively.
% Compared to the CI that uses zero as a lower bound for $\phi^2(\tau)$, the CIs widths have reductions of $0.04$, $0.09$ and $1.62$, respectively, when applying \cite{Aronow2013ACoUEotATEiRE}'s \cite{Ding2018DTEV}'s and our bounds. 
It can be seen that the reductions for the three methods are not very large comparing with the naive method. The reason may be that $N\hat{\phi}_{1}^{2}/n_{1} + N\hat{\phi}_{0}^{2}/n_{0}$ is too large compared to the estimator of the lower bound for $\phi^{2}(\tau)$ in this specific problem, where $\hat{\phi}_{1}^{2}$ and $\hat{\phi}_{0}^{2}$ are the estimators for $\phi^{2}(y_{1})$ and $\phi^{2}(y_{0})$, respectively. Notice that the CI width is proportional to $\sqrt{N\hat{\phi}_{1}^{2}/n_{1} + N\hat {\phi}_{0}^{2}/n_{0} - \hat{\phi}_{\rm B}^{2}}$, where $\hat{\phi}_{\rm B}^{2}$ is the estimator of the lower bound for $\phi^{2}(\tau)$. This implies that the lower bound estimator $ \hat{\phi}_{\rm B}^{2}$ does not play a important role in the CI width if $N\hat{\phi}_{1}^{2}/n_{1} + N\hat{\phi}_{0}^{2}/n_{0}$ is much larger than $\hat{\phi}^{2}_{\rm B}$.

%Although all these reductions are not that significant considering the magnitude of the widths, our method makes a greater reduction in the aspect of the width of the $95\%$ CI compared to existing methods.

\subsection{Application to JOBS II}
In this section, we apply our approach proposed in Section \ref{sec:noncom} to a  dataset from the randomized trial JOBS II \citep{vinokur1995}  for illustration. Data used in this subsection are available from https://www.icpsr.umich.edu/ web/ICPSR/studies/2739. The JOBS II intervention trial studied the efficacy of a job training intervention in preventing depression caused by job loss and in prompting high-quality reemployment. The treatment consisted of five half-day training seminars that enhance the participants' job search strategies. The control group receives a booklet with some brief tips.  After some screening procedures, 1801 respondents were enrolled in this study, with 552 and 1249 respondents in the control and treatment groups, respectively. Of the respondents assigned to the treatment group, only $54\%$ participated in the treatment. Thus there is a large proportion of noncompliance in this study. The parameter of interest is the LATE of the treatment on the depression score (larger score indicating severer depression). We use the gender, the initial risk status and the economic hardship as the covariates and apply our method proposed in Section \ref{sec:noncom}. The estimates for $\tilde{\phi}^{2}_{\rm L}$ and $\tilde{\phi}^2_{\rm H}$ are $0.23$ and $0.81$, respectively. The $95\%$ CIs constructed using the naive bound zero and $\tilde{\phi}_{\rm L}$ are $[-0.2428, 0.0271]$ and $[-0.2360, 0.0202]$, respectively. Our method shortens the CI by $0.014$. When testing the null hypothesis that $\theta_{\rm c} = 0$ against the alternative hypothesis that $\theta_{\rm c} < 0$, the naive method gives the $\rm p$-value of $0.059$ while our method gives the $\rm p$-value of $0.049$. Thus our method is able to detect the treatment effect at $0.05$ significance level while the naive method is not.

\section{Discussion}\label{sec:discuss}
In this paper, we establish sharp variance bounds for the widely-used difference-in-means estimator and Wald estimator in the presence of covariates in completely randomized experiments. These bounds can help to improve the performance of inference procedures based on normal approximation.
We do not impose any assumption on the support of outcomes, hence our results are general and are applicable to both binary and continuous outcomes. 
Variances of the difference-in-means estimator in matched pair randomized experiments  \citep{imai2008variance} and the Horvitz-Thompson estimator in stratified randomized experiments and clustered randomized experiments  \citep{miratrix2013adjusting,Mukerjee2018USTfFPStICIfCE,Middleton2015UEotATEiCRE} share similar unidentifiable term as our
consideration. Moreover, the unidentifiable phenomenon also appears in the asymptotic variance
of regression adjustment estimators, see \cite{Lin2013ANoRAtEDRFC}, \cite{Freedman2008ORAiEwST}, \cite{Bloniarz2016LAoTEEiRE}. The
insights in this paper are also applicable in these settings and we do not discuss this in detail to avoid cumbersome notations. In addition, it is of great interest to extend our work to randomized experiments with other randomization schemes such as $2^2$ factorial design \citep{Lu2017SRBCIf2FDwBO} or some other complex assignment mechanisms \citep{Mukerjee2018USTfFPStICIfCE}.
%%%%%%%%%%%%%%%%%%%%%%%%%%%%%%%%%%%%%%%%%%%%%%%%%%%%%%%%%%%%%%%%%%%%%%%%%%%%%%%%%%%%%%%%%%%%%%%%%%%%%%%%%%%%%%%%%%%%%%%%%%%%

\section*{Acknowledgements}
This research was supported by the National Natural Science Foundation of China (General
project 11871460 and project for Innovative Research Group 61621003), and a
grant from the Key Lab of Random Complex Structure and Data Science, CAS.
\bigskip

%%%%%%%%%%%%%%%%%%%%%%%%%%%%%%%%%%%%%%%%%%%%%%%%%%%%%%%%%%%%%%%%%%%%%%%%%%%%%%%%%%%%%%%%%%
\appendix
	The Appendix is organized as follows. In Section S1, we provide
the proof of Theorem  \ref{thm:cov}. In Section S2, we prove the relationship (2.2) in Section 2 in the main text. We prove Theorem \ref{thm:consistency cov} in Section S3. The proofs of Theorem \ref{thm: asy level cov} , \ref{thm:asymptotic dist compliance}, \ref{thm:compliance} and \ref{thm:consistency compliance} are provided in Section S4, S5, S6 and S7, respectively. Proof of Theorem \ref{thm: asy level compliance} is similar to that of Theorem \ref{thm: asy level cov} and hence is omitted. Section S8 contains further simulation results on the performance of the confidence intervals when the proposed lower bounds are attained.
\par

\setcounter{section}{0}
\setcounter{equation}{0}
\def\theequation{S\arabic{section}.\arabic{equation}}
\def\thesection{S\arabic{section}}
\def\thetable{S\arabic{table}}

\fontsize{12}{14pt plus.8pt minus .6pt}\selectfont

\section{Proof of Theorem 1}
\begin{proof}
	Note that 
	\[\phi^2(\tau) =\sum_{i=1}^{K}\pi_k\frac{1}{N_k}\sum_{w_i=\xi_k}(y_{1i}-y_{0i})^2-\mu^2(\tau),\]
	where $N_k=N\pi_{k}$.
	Letting $a_k(s)=F_{1\mid k}^{-1}(s/N_k)$ and $b_k(s)=F_{0\mid k}^{-1}(s/N_k)$ for $k=1,\dots,K$ and $s=1,\dots,N_k$, then we have
	\[\frac{1}{N_k}\sum_{w_i=\xi_k}(y_{1i}-y_{0i})^2=\frac{1}{N_k}\sum_{i=1}^{N_k}a_k(i)^2+\frac{1}{N_k}\sum_{i=1}^{N_k}b_k(i)^2-\frac{2}{N_k}\sum_{i=1}^{N_k}a_k(i)b_k(\Pi_k(i))\]
	where $\Pi_k$ is a permutation on $\{1,\dots,N_k\}$. By the rearrangement inequality, we have
	\[\sum_{i=1}^{N_k}a_k(i)b_k(N_k-k+1)\leq \sum_{i=1}^{N_k}a_k(i)b_k(\Pi_k(i))\leq \sum_{i=1}^{N_k}a_k(i)b_k(i).\]
	Thus
	\begin{align*}
		&\int_{0}^{1}(F_{1\mid k}^{-1}(u)-F_{0\mid k}^{-1}(u))^2du=\frac{1}{N_k}\sum_{i=1}^{N_k}(a_k(i)-b_k(i))^2\leq \frac{1}{N_k}\sum_{w_i=\xi_k}(y_{1i}-y_{0i})^2
		\\&\leq \frac{1}{N_k}\sum_{i=1}^{N_k}(a_k(i)-b_k(N_k+1-i))^2=\int_{0}^{1}(F_{1\mid k}^{-1}(u)-F_{0\mid k}^{-1}(1-u))^2du
	\end{align*}
	This proves the bound in Theorem \ref{thm:cov}.
	
	Next, we prove the sharpness of the bound. Let  $\mathcal{U}=\{\mathbf{U}^*=(y_1^*,y_0^*,w^*): P(w^*=\xi_k)=\pi_k\ \text{and}\ P(y_t^* \leq y\mid w^*=\xi_k)=F_{t\mid k}(y)\ \text{for}\ t=0,1\ \text{and} \ k=1,\dots,K\}$ be the set of all  populations whose covariate distribution and distributions of each potential outcome conditional on the covariate are all identical to those of $\mathbf{U}$. 
	We first prove that the established bound can be attained by some population among $\mathcal{U}$.
	
	For $k=1,\dots,K$, let $i_{k(1)}< \cdots<i_{k(N_k)}$ be the indices in $\mathcal{I}_k = \{i:w_{i} = \xi_k\}$ in increasing order. Define the population $\mathbf{U}^{\rm L}$ consisting of $N$ units with two potential outcomes $y_{1i}^{\rm L}$ and $y_{0i}^{\rm L}$ and a vector of covariates $w_i^{\rm L}$ associated with unit $i$ for $i=1,\dots,N$. Let $y_{1i_{k(j)}}^{\rm L}=a_k(j)$, $y_{0i_{k(j)}}^{\rm L}=b_k(j)$  and $w_{i_{k(j)}}^{\rm L}=\xi_k$ for $k=1\dots,K$ and $j=1,\dots,N_k$. Let $\tau_i^{\rm L}=y_{1i}^{\rm L}-y_{0i}^{\rm L}$ for $i=1,\dots,N$. Then $P(y_t^{\rm L}\leq y\mid w^{\rm L}=\xi_k)=F_{t\mid k}(y)$ and $P(w^{\rm L}=\xi_k)=\pi_k$ for $t=0,1$, $k=1,\dots,K$ and $y\in \mathbf{R}$ and $1/N_k\sum_{i\in \mathcal{I}_k}(y_{1i}^{\rm L}-y_{0i}^{\rm L})^2=1/N_k\sum_{j=1}^{N_k}(a_k(j)-b_k(j))^2=\int_{0}^{1}(F_{1\mid k}^{-1}(u)-F_{0\mid k}^{-1}(u))^2du$ for $k=1,\dots,K$. Thus $\phi^2(\tau^{\rm L})=\phi^2_{\rm L}$, which attains the lower bound for $\phi^2(\tau)$.
	
	Define $\mathbf{U}^{\rm H}$ similarly with $y_{1i_{k(j)}}^{\rm H}=a_k(j)$, $y_{1i_{k(j)}}^{\rm H}=b_k(N_k+1-j)$ and $w_{i_{k(j)}}^{\rm H}=\xi_k$ for $k=1\dots,K$ and $j=1,\dots,N_k$. Let $\tau_i^{\rm H}=y_{1i}^{\rm H}-y_{0i}^{\rm H}$, then $P(y_t^{\rm H}\leq y\mid w^{\rm H}=\xi_k)=F_{t\mid k}(y)$ and $P(w^{\rm H}=\xi_k)=\pi_k$ for $t=0,1$, $k=1,\dots,K$ and $y\in \mathbf{R}$. Moreover,  $1/N_k\sum_{i\in \mathcal{I}_k}(y_{1i}^{\rm H}-y_{0i}^{\rm H})^2=1/N_k\sum_{j=1}^{N_k}(a_k(j)-b_k(N_k+1-j))^2=\int_{0}^{1}(F_{1\mid k}^{-1}(u)-F_{0\mid k}^{-1}(1-u))^2du$ for $k=1,\dots,K$. Thus $\phi^2(\tau^{\rm H})=\phi^2_{\rm H}$, which attains the upper bound for $\phi^2(\tau)$.
	
	Next, we show the sharpness of the bound. We show the result for the lower bound, i.e., $\phi_{\rm L}^{2}$ is no smaller than any bound in $\mathcal{B}_{\rm L}$, and the result for the upper bound follows similarly.
	For any bound $b_{\rm L}$ in $\mathcal{B}_{\rm L}$, we have $b_{\rm L} = f(\pi_{1},\dots,\pi_{K}, F_{0\mid 1}, \dots, F_{0\mid K}, F_{1\mid 1}, \dots, F_{1\mid K})$ for some functional $f$ according to the definition of $\mathcal{B}_{\rm L}$. For any $\mathbf{U}^{*}$ in $\mathcal{U}$, define $\pi_{k}^{*}=P(w^{*}=\xi_k)$ and $F_{t\mid k}^{*}(y)=P(y_{t}^{*}\leq y\mid w=\xi_k)$ for $t=0,1$ and $k=1,\dots,K$.
	Let $\tau^* = (y_{11}^* - y_{01}^*, \dots, y_{1N}^* - y_{0N}^*)$. By the applying the bound to the population $\mathbf{U}^*$, we have 
	\begin{equation}\label{eq: bound on U*}
		f(\pi_{1}^{*},\dots,\pi_{K}^{*}, F_{0\mid 1}^{*}, \dots, F_{0\mid K}^{*}, F_{1\mid 1}^{*}, \dots, F_{1\mid K}^{*}) \leq \phi^{2}(\tau^{*}).
	\end{equation}
	Because $\mathbf{U}^{*} \in \mathcal{U}$, it holds that $\pi_{k}^{*}=\pi_{k}$ and $F_{t\mid k}^{*}(y)=F_{t\mid k}(y)$
	for $t=0,1$. Thus \eqref{eq: bound on U*} implies that 
	\begin{equation}\label{eq: bound order}
		\begin{aligned}
			b_{\rm L} &=  f(\pi_{1},\dots,\pi_{K}, F_{0\mid 1}, \dots, F_{0\mid K}, F_{1\mid 1}, \dots, F_{1\mid K})\\
			& = f(\pi_{1}^{*},\dots,\pi_{K}^{*}, F_{0\mid 1}^{*}, \dots, F_{0\mid K}^{*}, F_{1\mid 1}^{*}, \dots, F_{1\mid K}^{*}) \leq \phi^{2}(\tau^{*})
		\end{aligned}
	\end{equation}
	for any $\mathbf{U}^{*} \in \mathcal{U}$. We have shown that there is some $\mathbf{U}^{*} \in \mathcal{U}$ such that $\phi_{\rm L}^{2} = \phi^{2}(\tau^{*})$. Hence \eqref{eq: bound order} implies $b_{\rm L} \leq \phi_{\rm L}^{2}$ for any $b_{\rm L} \in \mathcal{B}_{\rm L}$, which completes the proof.
\end{proof}

\section{Proof of Theorem 2}\label{sec:proof of thm 2}
\begin{proof} Throughout this and the following proofs, for any real number $a$ and $b$, we let $a\land b=\min\{a,b\}$, $a\lor b=\min\{a,b\}$, $a_+ = a\land 0$ and $a_-=-(a\land 0)$. For any function $H$ and any constant $C$, we let
	\[H_{C}(y)=\left\{
	\begin{array}{lc}
		0&\ y<-C\\
		H(y)&\ -C\leq y<C\\
		1&\ y\geq C
	\end{array}
	\right.\]
	and $H_{C}^{-1}(y) = (H_{C})^{-1}(y)$.
	We prove consistency of the lower bound estimator only, 
	and the consistency  for the upper bound estimator follows similarly.
	For any $0 < \epsilon < 1$, it is easy to verify that $|\hat{\theta} - \theta| \leq \min\{1/(4\theta), 1/\sqrt{2}\}\epsilon$ implies $|\hat{\theta}^{2} - \theta^{2}| \leq \epsilon$. Thus $\mathbb{P}(|\hat{\theta}^{2}-\theta^{2}| > \epsilon) \leq \mathbb{P}(|\hat{\theta} - \theta| > \min\{1/(4\theta), 1/\sqrt{2}\}\epsilon) \leq \var(\hat{\theta}) \max\{16\theta^{2}, 2\}/\epsilon$ by Chebyshev's inequality. 
	Under Condition \ref{cond:fourth moment}, we have $\theta \leq (2C_{\rm M})^{1/4}$ and $\sigma^{2} \leq C_{\rm M}^{1/2}(1/n_{1} + 1/n_{0})N/(N-1)$ by Jensen's inequality. Thus for any $0<\epsilon, \delta<1$, we have 
	\begin{equation}\label{eq: hat_theta^2}
		\mathbb{P}(|\hat{\theta}^{2}-\theta^{2}| \leq \epsilon) \geq 1-\delta
	\end{equation}
	for sufficiently large $N$.
	Let $\hat{\pi}_{tk}=\sum_{T_{i} = t}1\{w=\xi_k\}/n_{t}$ for $t=0,1$ and $k=1,\dots,K$. Then for any $C > 0$, 
	\begin{align*}
		&\vert\sum_{k=1}^K\hat{\pi}_k\int_{0}^{1}(\hat{F}_{1\mid k}^{-1}(u)-\hat{F}_{0\mid k}^{-1}(u))^2du-\sum_{k=1}^{K}\pi_{k}\int_{0}^{1}(F_{1\mid k}^{-1}(u)-F_{0\mid k}^{-1}(u))^2du \vert& \\
		&\leq \vert\sum_{k=1}^K\hat{\pi}_k\int_{0}^{1}(\hat{F}_{1\mid k}^{-1}(u)-\hat{F}_{0\mid k}^{-1}(u))^2du-\sum_{k=1}^K\hat{\pi}_k\int_{0}^{1}(\hat{F}_{1\mid k,C}^{-1}(u)-\hat{F}_{0\mid k,C}^{-1}(u))^2du\vert&\\
		& + \vert\sum_{k=1}^{K}\pi_k\int_{0}^{1}(F_{1\mid k}^{-1}(u)-F_{0\mid k}^{-1}(u))^2du- \sum_{k=1}^{K}\pi_k\int_{0}^{1}(F_{1\mid k,C}^{-1}(u)-F_{0\mid k,C}^{-1}(u))^2du \vert&\\
		& + \vert\sum_{k=1}^K\hat{\pi}_k\int_{0}^{1}(\hat{F}_{1\mid k,C}^{-1}(u)-\hat{F}_{0\mid k,C}^{-1}(u))^2du-\sum_{k=1}^{K}\pi_k\int_{0}^{1}(F_{1\mid k,C}^{-1}(u)-F_{0\mid k,C}^{-1}(u))^2du \vert&\\
		&\equalscolon I_1+I_2+I_3.&
	\end{align*}
	Hence to prove Theorem \ref{thm:consistency cov}, it suffices to show $I_{1} + I_{2} + I_{3} \stackrel{P}{\to} 0$.
	Note that
	\begin{align*}
		I_2 &\leq  \vert\sum_{k=1}^{K}\pi_k\int_{0}^{1}(F_{1\mid k}^{-1}(u)-F_{1\mid k,C}^{-1}(u))(F_{1\mid k}^{-1}(u)+F_{1\mid k,C}^{-1}(u)-F_{0\mid k}^{-1}(u)-F_{0\mid k,C}^{-1}(u))du\vert\\
		& + \vert\sum_{k=1}^{K}\pi_k\int_{0}^{1}(F_{0\mid k}^{-1}(u)-F_{0\mid k,C}^{-1}(u))(F_{1\mid k}^{-1}(u)+F_{1\mid k,C}^{-1}(u)-F_{0\mid k}^{-1}(u)-F_{0\mid k,C}^{-1}(u))du\vert\\
		&\leq\Big(\sum_{k=1}^{K}\pi_k\int_{0}^{1}(F_{1\mid k}^{-1}(u)-F_{1\mid k,C}^{-1}(u))^2du\Big)^\frac{1}{2}\\
		&\times\Big(\sum_{k=1}^{K}\pi_k\int_{0}^{1}(F_{1\mid k}^{-1}(u)+F_{1\mid k,C}^{-1}(u)-F_{0\mid k}^{-1}(u)-F_{0\mid k,C}^{-1}(u))^2du\Big)^\frac{1}{2} \\
		& + \Big(\sum_{k=1}^{K}\pi_k\int_{0}^{1}(F_{0\mid k}^{-1}(u)-F_{0\mid k,C}^{-1}(u))^2du\Big)^\frac{1}{2}\\
		&\times\Big(\sum_{k=1}^{K}\pi_k\int_{0}^{1}(F_{1\mid k}^{-1}(u)+F_{1\mid k,C}^{-1}(u)-F_{0\mid k}^{-1}(u)-F_{0\mid k,C}^{-1}(u))^2du\Big)^\frac{1}{2},
	\end{align*}
	where the second inequality follows from Cauchy-Schwartz inequality.
	
	Under Condition \ref{cond:fourth moment}, we have $I_2\leq 8C_M^{3/4}/C$ because  \[
	\sum_{k=1}^{K}\pi_k\int_{0}^{1}(F_{1\mid k}^{-1}(u)-F_{1\mid k,C}^{-1}(u))^2du\leq \frac{1}{N}\sum_{\vert y_{ti}\vert\geq C}y_{ti}^2 \leq \frac{C_M}{C^2}\]
	for $t=0,1$ and
	\begin{align*}
		&\sum_{k=1}^{K}\pi_k\int_{0}^{1}(F_{1\mid k}^{-1}(u)+F_{1\mid k,C}^{-1}(u)-F_{0\mid k}^{-1}(u)-F_{0\mid k,C}^{-1}(u))^2du\\
		&\leq 4\sum_{k=1}^{K}\pi_k\int_{0}^{1}[(F_{1\mid k}^{-1}(u))^2+(F_{1\mid k,C}^{-1}(u))^2+(F_{0\mid k}^{-1}(u))^2+(F_{0\mid k,C}^{-1}(u))^2]du\\
		&\leq \frac{8}{N}\sum_{i=1}^{N}y_{1i}^2 + \frac{8}{N}\sum_{i=1}^{N}y_{0i}^2\leq 16\sqrt{C_M}.
	\end{align*}
	
	By similar arguments, we have 
	\begin{align*}
		I_1&\leq \Big(\sum_{k=1}^{K}\hat{\pi}_k\int_{0}^{1}(\hat{F}_{1\mid k}^{-1}(u)-\hat{F}_{1\mid k,C}^{-1}(u))^2du\Big)^\frac{1}{2}\\
		&\times\Big(4\sum_{k=1}^{K}\hat{\pi}_k\int_{0}^{1}[(\hat{F}_{1\mid k}^{-1}(u))^2+(\hat{F}_{1\mid k,C}^{-1}(u))^2+(\hat{F}_{0\mid k}^{-1}(u))^2+(\hat{F}_{0\mid k,C}^{-1}(u))^2]du\Big)^\frac{1}{2} \\
		&+ \Big(\sum_{k=1}^{K}\hat{\pi}_k\int_{0}^{1}(\hat{F}_{0\mid k}^{-1}(u)-\hat{F}_{0\mid k,C}^{-1}(u))^2du\Big)^\frac{1}{2}\\
		&\times\Big(4\sum_{k=1}^{K}\hat{\pi}_k\int_{0}^{1}[(\hat{F}_{1\mid k}^{-1}(u))^2+(\hat{F}_{1\mid k,C}^{-1}(u))^2+(\hat{F}_{0\mid k}^{-1}(u))^2+(\hat{F}_{0\mid k,C}^{-1}(u))^2]du\Big)^\frac{1}{2}.
	\end{align*}
	Because for $t=0,1$, \[\sum_{k=1}^{K}\hat{\pi}_k\int_{0}^{1}(\hat{F}_{t\mid k}^{-1}(u)-\hat{F}_{t\mid k,C}^{-1}(u))^2du\leq \max_k\frac{\hat{\pi}_k}{\hat{\pi}_{tk}}\frac{N}{n_t}\frac{1}{N}\sum_{\vert y_{ti}\vert \geq C}y_{ti}^2
	\leq \frac{N}{n_t} \max_k\frac{\hat{\pi}_k}{\hat{\pi}_{tk}}\frac{C_M}{C^2},\]
	and
	\begin{align*}
		& \sum_{k=1}^{K}\hat{\pi}_k\int_{0}^{1}[(\hat{F}_{1\mid k}^{-1}(u))^2+(\hat{F}_{1\mid k,C}^{-1}(u))^2+(\hat{F}_{0\mid k}^{-1}(u))^2+(\hat{F}_{0\mid k,C}^{-1}(u))^2]du &\\
		&\leq 
		\frac{N}{n_1\land n_0}\max_{t,k}\frac{\hat{\pi}_k}{\hat{\pi}_{tk}}\Big(\frac{2}{N}\sum_{i=1}^{N}y_{1i}^2
		+ \frac{2}{N}\sum_{i=1}^{N}y_{0i}^2\Big)&  \\
		&\leq \frac{4N}{n_1\land n_0}\max_{t,k}\frac{\hat{\pi}_k}{\hat{\pi}_{tk}}\sqrt{C_M}, &
	\end{align*}
	we have
	\[
	I_1\leq \frac{8N}{n_1\land n_0}\max_{t,k}\frac{\hat{\pi}_k}{\hat{\pi}_{tk}}\frac{C_M^{\frac{3}{4}}}{C}.
	\]
	For the last term $I_3$, we have
	\begin{align*}
		I_3 & =  \left\vert\sum_{k=1}^K\hat{\pi}_k\int_{0}^{1}(\hat{F}_{1\mid k,C}^{-1}(u)-\hat{F}_{0\mid k,C}^{-1}(u))^2du-\sum_{k=1}^{K}\pi_k\int_{0}^{1}(F_{1\mid k,C}^{-1}(u)-F_{0\mid k,C}^{-1}(u))^2du \right\vert\\
		&\leq  \left\vert\sum_{k=1}^{K}\hat{\pi}_k\Big(\int_{0}^{1}(\hat{F}_{1\mid k,C}^{-1}(u)-\hat{F}_{0\mid k,C}^{-1}(u))^2du -  \int_{0}^{1}(F_{1\mid k,C}^{-1}(u)-F_{0\mid k,C}^{-1}(u))^2du\Big) \right\vert\\
		&+ \left\vert\sum_{k=1}^{K}(\hat{\pi}_k-\pi_k)\int_{0}^{1}(F_{1\mid k,C}^{-1}(u)-F_{0\mid k,C}^{-1}(u))^2du\right\vert\\
		&\leq \max_k\left\vert\int_{0}^{1}(\hat{F}_{1\mid k,C}^{-1}(u)-\hat{F}_{0\mid k,C}^{-1}(u))^2du -  \int_{0}^{1}(F_{1\mid k,C}^{-1}(u)-F_{0\mid k,C}^{-1}(u))^2du\right\vert\sum_{k=1}^{K}\hat{\pi}_{k} \\
		&+ \left\vert\sum_{k=1}^{K}(\frac{\hat{\pi}_k}{\pi_{k}}-1)\pi_{k}\int_{0}^{1}(F_{1\mid k,C}^{-1}(u)-F_{0\mid k,C}^{-1}(u))^2du\right\vert
		\\
		&\leq  \max_k\left\vert\int_{0}^{1}(\hat{F}_{1\mid k,C}^{-1}(u)-\hat{F}_{0\mid k,C}^{-1}(u))^2du -  \int_{0}^{1}(F_{1\mid k,C}^{-1}(u)-F_{0\mid k,C}^{-1}(u))^2du\right\vert\\
		&+ \max_k\left\vert\frac{\hat{\pi}_{k}}{\pi_k}-1\right\vert \left\vert\sum_{k=1}^{K}\pi_k\int_{0}^{1}(F_{1\mid k,C}^{-1}(u)-F_{0\mid k,C}^{-1}(u))^2du\right\vert\\
		&\equalscolon I_{31} + I_{32}.
	\end{align*}
	By Condition \ref{cond:fourth moment}, we have
	\begin{align*}
		&\sum_{k=1}^{K}\pi_k\int_{0}^{1}(F_{1\mid k,C}^{-1}(u)-F_{0\mid k,C}^{-1}(u))^2du\\
		&\leq 
		2\sum_{k-1}^{K}\pi_k\int_0^1(F^{-1}_{1\mid k}(u))^2du + 2\sum_{k-1}^{K}\pi_k\int_0^1(F^{-1}_{0\mid k}(u))^2du \\
		& = 2\Big(\frac{1}{N}\sum_{i=1}^{N}y_{1i}^2+\frac{1}{N}\sum_{i=1}^{N}y_{0i}^2\Big)\\
		&\leq 2\sqrt{C_M}.
	\end{align*}
	Thus
	\[I_{32} \leq 2 \max_k \left\vert\frac{\hat{\pi}_k}{\pi_k}-1\right\vert \sqrt{C_M}.\]
	
	For $k=1,\dots,K$, $\int_{0}^{1}(F_{1\mid k,C}^{-1}(u)-F_{0\mid k,C}^{-1}(u))^2du$ is the square of the Wasserstein distance induced by $L_2$ norm between $F_{1\mid k,C}$ and
	$F_{0\mid k,C}$. 
	By the representation theorem \citep{bobkov2019one}[Theorem 2.11],
	\begin{align*}
		&\int_{0}^{1}(\hat{F}_{1\mid k,C}^{-1}(u)-\hat{F}_{0\mid k,C}^{-1}(u))^2du\\
		&=2\iint_{v\leq w}[(\hat{F}_{1\mid k,C}(v)-\hat{F}_{0 \mid k, C}(w))_+ + 
		(\hat{F}_{0\mid k,C}(v)-\hat{F}_{1 \mid k, C}(w))_+]dvdw.
	\end{align*}
	Because  $\hat{F}_{1\mid k,C}^{-1}(v)-\hat{F}_{0\mid k,C}^{-1}(w)\leq 0$ and
	$\hat{F}_{0\mid k,C}(v)-\hat{F}_{1 \mid k, C}(w)\leq 0$ when $v < -C$ or $w \geq C$, 
	the integral domain can be restricted to $\{-C\leq v \leq w < C\}$ without changing
	the integral.
	
	Similarly, we have 
	\begin{align*}
		&\int_{0}^{1}(F_{1\mid k,C}^{-1}(u)-F_{0\mid k,C}^{-1}(u))^2du\\
		&=2\iint_{-C\leq v\leq w < C}[(F_{1\mid k,C}(v)-F_{0 \mid k, C}(w))_+ + 
		(F_{0\mid k,C}(v)-F_{1 \mid k, C}(w))_+]dvdw.
	\end{align*}
	Because $\vert(u_1)_{+} - (u_2)_+\vert \leq \vert u_1 - u_2\vert$ for any $u_1$, $u_2$,
	\begin{align*}
		I_{31} &\leq  2\max_k\iint_{-C\leq v\leq w < C}[\vert\hat{F}_{1\mid k,C}(v) - F_{1\mid k,C}(v)\vert + \vert\hat{F}_{0\mid k,C}(v) - F_{0\mid k,C}(v)\vert
		\\
		&  \phantom{\leq2\max_k\iint_{-C\leq v\leq w < C}} + \vert\hat{F}_{1\mid k,C}(w) - F_{1\mid k,C}(w)\vert + 
		\vert\hat{F}_{0\mid k,C}(w) - F_{0\mid k,C}(w)\vert] dvdw \\
		&\leq 2C^2\max_k\{\sup_v\vert\hat{F}_{1\mid k,C}(v) - F_{1\mid k,C}(v)\vert + 
		\sup_v\vert\hat{F}_{0\mid k,C}(v) - F_{0\mid k,C}(v)\vert\}\\
		&\leq 4C^2\max_{t,k}\sup_v\vert\hat{F}_{t\mid k,C}(v) - F_{t\mid k,C}(v)\vert.
	\end{align*}
	
	For any given $M$ and $t=0,1$, let $s_{tk,j}=F_{t\mid k,C}^{-1}(j/M)$. By the standard technique in the proof of Glivenko-Cantelli theorem \citep{Vaart2000AS},
	\[\sup_v\vert\hat{F}_{t\mid k,C}(v) - F_{t\mid k,C}(v)\vert\leq \max_j\vert\hat{F}_{t\mid k,C}(s_{tk,j}) - F_{t\mid k,C}(s_{tk,j})\vert + \frac{1}{M}.\]
	For $t=0,1$ and $i=1,\dots,N$, let $y_{ti}^*=y_{ti}1\{-C\leq y_{ti} < C\} + C1\{y_{ti}\geq C\} - C1\{y_{ti} < -C\}$.
	Then
	\begin{align*}
		&\hat{F}_{tk,C}(s_{tk,j})-F_{tk,C}(s_{tk,j})\\
		& = \frac{1}{\hat{\pi}_{tk}}\frac{1}{n_t}\sum_{T_i=t}1\{y_{ti}^*\leq s_{tk,j},x=\xi_k\}-\frac{1}{\pi_k}\frac{1}{N}\sum_{i=1}^{N}1\{y_{ti}^*\leq s_{tk,j},x=\xi_k\}\\
		& =  \frac{1}{\pi_k}\left(\frac{1}{n_t}\sum_{T_i=t}1\{y_{ti}^*\leq s_{tk,j},x=\xi_k\} - \frac{1}{N}\sum_{i=1}^{N}1\{y_{ti}^*\leq s_{tk,j},x=\xi_k\}\right)\\
		& -
		\frac{1}{\pi_k\hat{\pi}_{tk}}(\hat{\pi}_{tk}-\pi_k)\frac{1}{n_t}\sum_{T_i=t}1\{y_{ti}^*\leq s_{tk,j},x=\xi_k\}.
	\end{align*}
	Because
	\[
	\left\vert\frac{1}{\hat{\pi}_{tk}n_t}\sum_{T_i=t}1\{y_{ti}^*\leq s_{tk,j},x=\xi_k\}\right\vert \leq 1,
	\]
	we have
	\begin{align*}
		I_{31} &\leq  4C^2\max_{t,k,j}\Big\{\Big\vert\frac{1}{\pi_k}\Big(\frac{1}{n_t}\sum_{T_i=t}1\{y_{ti}^*\leq s_{tk,j},x=\xi_k\}\\
		&\quad - \frac{1}{N}\sum_{i=1}^{N}1\{y_{ti}^*\leq s_{tk,j},x=\xi_k\}\Big)\Big\vert + 
		\frac{1}{\pi_k}\vert\hat{\pi}_{tk}-\pi_k\vert \Big\}
		+ \frac{1}{M}.
	\end{align*}
	
	For any small positive number $\epsilon$, one can choose $C$ and $M$ such that $8C_M^{3/4}/C\leq \epsilon$ and $1/M\leq \epsilon$. Then by Hoeffding inequality for sample without replacement \citep{bardenet2015} and the Bonferroni inequality, we have 
	\begin{align*}
		&\mathbb{P}\left(4C^2\max_{t,k,j}\Big\{\frac{1}{\pi_k}\Big\vert\frac{1}{n_t}\sum_{T_i=t}1\{y_{ti}^*\leq s_{tk,j},x=\xi_k\} - \frac{1}{N}\sum_{i=1}^{N}1\{y_{ti}^*\leq s_{tk,j},x=\xi_k\}\Big\vert\Big\}\geq \epsilon\right)\\
		& \leq 4MK\exp\left(-\frac{n_1\land n_0\epsilon^2\min_k\pi_k^2}{8C^4}\right),
	\end{align*}
	
	and
	\[\mathbb{P}\left( \max_{t,k}\frac{1}{\pi_k}\vert\hat{\pi}_{tk}-\pi_k\vert\geq \epsilon\right) \leq 2K\exp\left( - 2 n_1\land n_0\epsilon^2\min_k\pi_k^2\right).\]
	By Conditions \ref{cond:positivity} and \ref{cond:order1} the right hand side of these two inequalities converge to zero because 
	\[n_1\land n_0\min_k\pi_k^2 -  C^{*}\log K\geq C_\pi (n_1\land n_0/N)(N/K^2) - C^{*} \log K \to \infty,\]
	with $C^{*} = 8C^{4} \lor (1/2)$.
	
	Hence,  for any $\delta > 0$ and sufficiently large $N$, $I_{31}\leq 3\epsilon$ and
	\[\max_{t,k}\left\vert\frac{\hat{\pi}_k}{\hat{\pi}_{tk}}\right\vert\leq \frac{1 + \epsilon}{1 - \epsilon}\]
	with probability greater than $1-\delta$.
	Without loss of generality, we let $\epsilon \leq 1/3$, then $(1 + \epsilon)/(1 - \epsilon) \leq 2$ and for sufficiently large $N$ we have
	\begin{align*}&\left\vert\sum_{k=1}^K\hat{\pi}_k\int_{0}^{1}(\hat{F}_{1\mid k}^{-1}(u)-\hat{F}_{0\mid k}^{-1}(u))^2du-\sum_{k=1}^{K}\int_{0}^{1}(F_{1\mid k}^{-1}(u)-F_{0\mid k}^{-1}(u))^2du \right\vert\\
		&\leq \left(\frac{4}{\rho_1\land\rho_0}+2\sqrt{C_M}+4\right)\epsilon
	\end{align*}
	with probability greater than $1-\delta$. Combining this with \eqref{eq: hat_theta^2} completes the proof of  Theorem \ref{thm:consistency cov}.
\end{proof}

\section{Proof of Theorem 3}
In this proof, we use $\epsilon,\delta$ to denote small positive numbers whose values may change from place to place. Recall the definition of $\mathcal{P}_{N}$ in \eqref{eq: def pop class} in the main text. We assume without loss of generality that $\mathcal{P}_{N}$ is non-empty.
Under Condition \ref{cond:order1}, according to the proof of Theorem \ref{thm:consistency cov}, for $\mathbf{U}_{N}\in \mathcal{P}_{N}$ and any $\epsilon, \delta>0$, we have $\mathbb{P}(|\hat{\phi}_{\rm L}^{2} - \phi_{\rm L}^{2}| \leq \epsilon) \geq 1-\delta$
for any $N$ larger than a threshold $T_{1}$. Note that in the proof of Theorem  \ref{thm:consistency cov}, the threshold $T_{1}$ can be chosen to be uniform in $\mathbf{U}_{N}\in \mathcal{P}_{N}$. That is, for any small positive numbers $\epsilon, \delta$, there is some $T_{1}$ such that for $N>T_{1}$ we have 
\begin{equation}\label{eq: uniform consistency}
	\inf_{\mathbf{U}_{N}\in \mathcal{P}_{N}}\mathbb{P}(|\hat{\phi}_{\rm L}^{2} - \phi_{\rm L}^{2}| \leq \epsilon) \geq 1-\delta.
\end{equation}
Recall the definition of $\hat{\sigma}^{2}$ in Section \ref{subsec: estimation cov} and the definition of $\hat{\phi}^{2}_{1}$, $\hat{\phi}^{2}_{0}$ in Section \ref{subsec:preliminaries} in the main text.
According to Theorem \ref{thm:cov},
\begin{equation}\label{eq: pointwise lower bound prob}
	\begin{aligned}
		\mathbb{P}(\hat{\sigma}^{2} \geq \sigma^{2} - \epsilon)& \geq 1 - \mathbb{P}\left(\frac{N}{n_{1}}|\phi^{2}(y_{1}) - \hat{\phi}^{2}_{1}| \geq \frac{\epsilon}{3}\right) - \mathbb{P}\left(\frac{N}{n_{0}}|\phi^{2}(y_{1}) - \hat{\phi}^{2}_{0}| \geq \frac{\epsilon}{3}\right)\\
		& \quad - \mathbb{P}\left(|\hat{\phi}_{\rm L}^{2} - \phi_{\rm L}^{2}| \geq \frac{\epsilon}{3}\right).
	\end{aligned}
\end{equation}
Thus for any  $\epsilon, \delta > 0$, there is some threshold $T_{2}$, for $N > T_{2}$, we have
\begin{equation}\label{eq: uniform lower bound prob}
	\inf_{\mathbf{U}_{N}\in \mathcal{P}_{N}} \mathbb{P}(\hat{\sigma}^{2} \geq \sigma^{2} - \epsilon) \geq 1-\delta
\end{equation}
by \eqref{eq: uniform consistency}, \eqref{eq: pointwise lower bound prob}, Chebyshev's inequality and straightforward calculation of mean and variance of $\hat{\phi}_{1}^{2}$, $\hat{\phi}_{0}^{2}$. According to the Cauchy-Schwartz inequality, we have
\begin{equation}\label{eq: bound sigma}
	\begin{aligned}
		\sigma^{2} 
		&= \frac{N}{n_1}\phi^2(y_1)+\frac{N}{n_0}\phi^2(y_0)-\phi^2(\tau)\\
		&\geq
		\frac{N}{n_1}\phi^2(y_1)+\frac{N}{n_0}\phi^2(y_0)-\phi^2(y_{1}) - \phi^2(y_{0}) \\
		&\geq L_{2}\frac{n_{0}}{n_{1}} + L_{2}\frac{n_{1}}{n_{0}}\\
		&> L_{2}
	\end{aligned}
\end{equation}
if $\mathbf{U}_{N} \in \mathcal{P}_{N}$. Then \eqref{eq: uniform lower bound prob}
and \eqref{eq: bound sigma} imply 
\begin{equation}\label{eq: uniform lower bound prob sqrt}
	\inf_{\mathbf{U}_{N}\in \mathcal{P}_{N}} \mathbb{P}(\hat{\sigma} \geq \sigma - \epsilon) \geq 1-\delta
\end{equation}
for any $N > T_{3}$, where $T_{3}$ is some threshold that depends on $\epsilon$, $\delta$.

By the finite population central limit theorem \citep{Freedman2008ORAiEwST}[Theorem 1], for any sequence of finite populations $\{\mathbf{U}_{N}\}_{N=2}^{\infty}$ such that $\mathbf{U}_{N} \in \mathcal{P}_{N}$, we have
\[\sqrt{N}\sigma^{-1}(\hat{\theta} - \theta) \stackrel{d}{\to} N(0,1)\]
as $N \to \infty$.  This and \eqref{eq: bound sigma} imply for any $0< \delta < \alpha < 1$ and any sequence of finite populations $\{\mathbf{U}_{N}\}_{N=2}^{\infty}$ such that $\mathbf{U}_{N} \in \mathcal{P}_{N}$, there is some $0 < \epsilon <  \sqrt{L_{2}}$ such that
\begin{equation}\label{eq: oracle pointwise coverage rate}
	\begin{aligned}
		&\mathop{\lim \inf}\limits_{N\to\infty}\mathbb{P}\left(\theta \in \left[\hat{\theta} - q_{\frac{\alpha}{2}}(\sigma - \epsilon) N^{-\frac{1}{2}}, \hat{\theta} + q_{\frac{\alpha}{2}}(\sigma - \epsilon) N^{-\frac{1}{2}}\right]\right) \\ 
		&= \mathop{\lim \inf}\limits_{N\to\infty}\mathbb{P}\left(\sqrt{N}\sigma^{-1}(\hat{\theta} - \theta) \in \left[- q_{\frac{\alpha}{2}}(1 - \sigma^{-1}\epsilon), q_{\frac{\alpha}{2}}(1 - \sigma^{-1}\epsilon)\right]\right)\\
		&\geq \alpha - \delta,
	\end{aligned}
\end{equation} 
where $q_{\alpha / 2}$ is the upper $\alpha/2$ quantile of a standard normal distribution.
Next we show that 
\begin{equation}\label{eq: oracle uniform coverage rate}
	\mathop{\lim \inf}\limits_{N\to\infty}\inf_{\mathbf{U}_{N} \in \mathcal{P}_{N}}\mathbb{P}\left(\theta \in \left[\hat{\theta} - q_{\frac{\alpha}{2}}(\sigma - \epsilon)N^{-\frac{1}{2}}, \hat{\theta} + q_{\frac{\alpha}{2}}(\sigma - \epsilon) N^{-\frac{1}{2}}\right]\right) \geq \alpha - \delta.
\end{equation}
We prove this by contradiction. If \eqref{eq: oracle uniform coverage rate} does not hold, then there is some $c_{*}>0$, for any $N_{0}>0$, there is some $N_{1} > N_{0}$ and $\mathbf{U}_{N_{1}}^{\prime}$ such that
\[\mathbb{P}\left(\theta \in \left[\hat{\theta} - q_{\frac{\alpha}{2}}(\sigma - \epsilon) N^{-\frac{1}{2}}, \hat{\theta} + q_{\frac{\alpha}{2}}(\sigma - \epsilon) N^{-\frac{1}{2}}\right]\right) \leq \alpha - \delta - c_{*}\]
under $\mathbf{U}_{N_{1}}^{\prime}$. Following the same procedure, we can extract a subsequence $\{\mathbf{U}_{N_{m}}^{\prime}\}_{m=1}^{\infty}$ with $N_{1} < N_{2} < \dots$ and $\mathbf{U}_{N_{m}}^{\prime} \in \mathcal{P}_{N_{m}}$ such that 
\[\mathbb{P}\left(\theta \in \left[\hat{\theta} - q_{\frac{\alpha}{2}}(\sigma - \epsilon) N^{-\frac{1}{2}}, \hat{\theta} + q_{\frac{\alpha}{2}}(\sigma - \epsilon) N^{-\frac{1}{2}}\right]\right) \leq \alpha - \delta - c_{*}\]
under $\mathbf{U}_{N_{m}}^{\prime}$ for $m=1,2,\dots$ Then it follows that for any sequence of finite populations $\{\mathbf{U}_{N}\}_{N=2}^{\infty}$ such that $\{\mathbf{U}_{N}\}_{N=2}^{\infty}$ contains $\{\mathbf{U}_{N_{m}}^{\prime}\}_{m=1}^{\infty}$ as a subsequence, we have 
\[	\lim\limits_{N\to\infty}\mathbb{P}\left(\theta \in \left[\hat{\theta} - q_{\frac{\alpha}{2}}(\sigma - \epsilon) N^{-\frac{1}{2}}, \hat{\theta} + q_{\frac{\alpha}{2}}(\sigma - \epsilon) N^{-\frac{1}{2}}\right]\right) \leq \alpha - \delta - c_{*},\]
which is contradict with \eqref{eq: oracle pointwise coverage rate}. This proves \eqref{eq: oracle uniform coverage rate}. According to \eqref{eq: oracle uniform coverage rate}, there is some threshold $T_{3}$ such that 
\[
\inf_{\mathcal{P}_{N}}\mathbb{P}\left(\theta \in \left[\hat{\theta} - q_{\frac{\alpha}{2}}(\sigma - \epsilon)N^{-\frac{1}{2}}, \hat{\theta} + q_{\frac{\alpha}{2}}(\sigma - \epsilon) N^{-\frac{1}{2}}\right]\right) \geq \alpha - 2\delta
\]
for $N > T_{3}$.
Combining this with \eqref{eq: uniform lower bound prob}, we have 
\[\inf_{\mathcal{P}_{N}}\mathbb{P}\left(\theta \in \left[\hat{\theta} - q_{\frac{\alpha}{2}}\hat{\sigma}N^{-\frac{1}{2}}, \hat{\theta} + q_{\frac{\alpha}{2}}\hat{\sigma} N^{-\frac{1}{2}}\right]\right) \geq \alpha - 3\delta\]
for $N > \max\{T_{2}, T_{3}\}$. This proves the theorem since $\delta$ is an arbitrary small positive number.

\section{Proof of Theorem 4}
\begin{proof}
	Under Conditions \ref{cond:fourth moment} and \ref{cond:nonsingular cov}, the conditions of the finite population central limit theorem \citep{Freedman2008ORAiEwST}[Theorem 1] is satisfied and hence
	\[\sqrt{N}V_{N}^{-1}\left(\frac{1}{n_1}\sum_{T_i=1}y_{1i} - \mu(y_1), \frac{1}{n_1}\sum_{T_i=1}d_{1i} - \mu(d_1), \frac{1}{n_0}\sum_{T_i=0}y_{0i} - \mu(y_0), \frac{1}{n_0}\sum_{T_i=0}d_{0i} - \mu(d_0)\right)^{\T}\]
	converges to a multivariate normal distribution, where $V_{N}$ is defined in Condition \ref{cond:nonsingular cov}. Note that we have the decomposition
	\begin{align*}
		\hat{\theta}_{\rm c} - \theta_{\rm c} =\hat{\pi}_{\rm c}^{-1}(\hat{\theta} - \hat{\pi}_{\rm c}\theta_{\rm c})
		=  \left(\hat{\pi}_{\rm c}^{-1} - \pi_{\rm c}^{-1}\right)\left(\hat{\theta} - \hat{\pi}_{\rm c}\theta_{\rm c}\right) + \pi_{\rm c}^{-1}(\hat{\theta} - \hat{\pi}_{\rm c}\theta).
	\end{align*}
	Under Condition \ref{cond:fourth moment}, it is easy to show that the trace and hence the spectral norm of $V_{N}$ is bounded.
	By the strong instrument assumption, the fact that $V_{N}$ has bounded spectral norm and the asymptotic normality invoked before, we have 
	\[
	\hat{\pi}_{\rm c}^{-1} - \pi_{\rm c}^{-1} = O_{p}\left(\frac{1}{\sqrt{N}}\right)
	\] 
	and 
	\[
	\hat{\theta} - \hat{\pi}_{\rm c}\theta_{\rm c} = O_{p}\left(\frac{1}{\sqrt{N}}\right).
	\] 
	Hence
	\[
	\hat{\theta}_{\rm c} - \theta_{\rm c} = \pi_{\rm c}^{-1}(\hat{\theta} - \hat{\pi}_{\rm c}\theta_{\rm c}) + o_p\left(\frac{1}{\sqrt{N}}\right).
	\]
	Straightforward calculation can show that $(\hat{\theta} - \hat{\pi}_{\rm c}\theta_{\rm c})/\pi_{\rm c}$ has mean zero and variance 
	\[\frac{1}{\pi_{\rm c}^{2}(N-1)}\left(\frac{N}{n_1}\phi^2(\tilde{y}_{1})+\frac{N}{n_0}\phi^2(\tilde{y}_{0})-\phi^2(\tilde{\tau})\right) = \frac{\sigma^{2}_{\rm c}}{N - 1}.\]
	Again by the finite population central limit theorem we have
	\[\sqrt{N}\sigma_{\rm c}^{-1}(\hat{\theta}_{\rm c} - \theta_{\rm c}) = \sqrt{N}\sigma_{\rm c}^{-1}\pi_{\rm c}^{-1}(\hat{\mu}(\tau) - \hat{\pi}_{\rm c}\theta_{\rm c}) + o_p(1) \overset{d}{\to}N(0, 1).\]
\end{proof}

\section{Proof of Theorem 5}
\begin{proof}
	According to exclusion restriction, if $g_i = \rm a$ or $\rm n$, we have 
	\[\tilde{\tau}_i = y_{1i} - y_{0i} - \theta_{\rm c}(d_{1i} - d_{0i}) = 0.\]
	By the definition of $\theta_{\rm c}$, we have $\mu(\tilde{\tau}) = 0$.
	Thus 
	\[\phi^2(\tilde{\tau}) = \frac{1}{N}\sum_{i=1}^{N} \tilde{\tau}^{2}_{i} = \pi_{\rm c} \frac{1}{N_{\rm c}}\sum_{g_i=\rm c}\tilde{\tau}_{i}^2\]
	where $N_{\rm c}$. Then the bound can be proved following the same arguments as in Theorem \ref{thm:cov}. Letting  $\mathcal{U}_{\rm c}=\{\mathbf{U}_{\rm c}^*=(y_1^*,y_0^*,w^*,g^*):
	\mathbf{U}_{\rm c}^* \text{ satisfies Assumption \ref{ass:iv}};$ \\
	$P(g^*=h) = \pi_{h},\ P(w^*=\xi_k\mid g^* = h)=\pi_{k\mid h}\ \text{and}\ P(y_t^* \leq y\mid w^*=\xi_k,g^*=h)=F_{t\mid (k,h)}(y)\
	\text{for}\ t=0,1\ \text{and} \ h={\rm a,c,n}\}$, the sharpness of the bound can be proved similarly as in Theorem \ref{thm:cov}.
\end{proof}

\section{Proof of Theorem 6} 
\begin{proof}
	Firstly, we provide some relationship that is useful in the proof.
	According to the monotonicity and exclusion restriction, we have 
	\[1\{g_i = {\rm c}\} = d_{1i} - d_{0i},\] \[(1 - d_{1i})1\{\tilde{y}_{0i}\leq y\} = (1 - d_{1i})1\{\tilde{y}_{1i}\leq y\}\] 
	and 
	\[d_{0i}1\{\tilde{y}_{1i}\leq y\} = d_{0i}1\{\tilde{y}_{0i}\leq y\}.\]
	Thus 
	\begin{align*}
		&\pi_{k\mid{\rm c}}\\
		& = \frac{\sum_{i=1}^{N}1\{g_i = {\rm c}\}1\{w_{i}=\xi_{k}\}}{\sum_{i=1}^{N}1\{g_{i}={\rm c}\}}\\
		& = \pi_{\rm c}^{-1} \left(\frac{1}{N}\sum_{i=1}^{N}d_{1i}1\{w_{i}=\xi_{k}\} - \frac{1}{N}\sum_{i=1}^{N}d_{0i}1\{w_{i}=\xi_{k}\}\right)\\
		& = \pi_{\rm c}^{-1} \left(\frac{1}{N}\sum_{i=1}^{N}(1 - d_{0i})1\{w_{i}=\xi_{k}\} - \frac{1}{N}\sum_{i=1}^{N}(1 - d_{1i})1\{w_{i}=\xi_{k}\}\right),
	\end{align*}
	\begin{align*}
		&\tilde{F}_{1\mid k}(y)\\
		& = \frac{\sum_{i=1}^{N}(d_{1i} - d_{0i})1\{\tilde{y}_{1i}\leq y\}1\{w_{i} = \xi_{k}\}}{\sum_{i=1}^{N}(d_{1i} - d_{0i})1\{w_{i} = \xi_{k}\}}\\
		& = \pi_{\rm c}^{-1}\pi_{k\mid{\rm c}}^{-1}\left(\frac{1}{N}\sum_{i=1}^{N}d_{1i}1\{\tilde{y}_{1i}\leq y\}1\{w_{i}=\xi_{k}\} - \frac{1}{N}\sum_{i=1}^{N}d_{0i}1\{\tilde{y}_{1i}\leq y\}1\{w_{i}=\xi_{k}\}\right)\\
		& = \pi_{\rm c}^{-1}\pi_{k\mid{\rm c}}^{-1}\left(\frac{1}{N}\sum_{i=1}^{N}d_{1i}1\{\tilde{y}_{1i}\leq y\}1\{w_{i}=\xi_{k}\} - \frac{1}{N}\sum_{i=1}^{N}d_{0i}1\{\tilde{y}_{0i}\leq y\}1\{w_{i}=\xi_{k}\}\right),
	\end{align*}
	and
	\begin{align*}
		&\tilde{F}_{0\mid k}(y)\\
		& = \frac{\sum_{i=1}^{N}(d_{1i} - d_{0i})1\{\tilde{y}_{0i}\leq y\}1\{w_{i} = \xi_{k}\}}{\sum_{i=1}^{N}(d_{1i} - d_{0i})1\{w_{i} = \xi_{k}\}}\\
		& =  \pi_{\rm c}^{-1}\pi_{k\mid{\rm c}}^{-1}\left(\frac{1}{N}\sum_{i=1}^{N}(1-d_{0i})1\{\tilde{y}_{0i}\leq y\}1\{w_{i}=\xi_{k}\}-\frac{1}{N}\sum_{i=1}^{N}(1 - d_{1i})1\{\tilde{y}_{0i}\leq y\}1\{w_{i}=\xi_{k}\}\right)\\
		& = \pi_{\rm c}^{-1}\pi_{k\mid{\rm c}}^{-1}\left(\frac{1}{N}\sum_{i=1}^{N}(1-d_{0i})1\{\tilde{y}_{0i}\leq y\}1\{w_{i}=\xi_{k}\} - \frac{1}{N}\sum_{i=1}^{N}(1 - d_{1i})1\{\tilde{y}_{1i}\leq y\}1\{w_{i}=\xi_{k}\}\right).
	\end{align*}
	Since the estimators in Theorem \ref{thm:consistency compliance} have the similar structure as those in Theorem \ref{thm:consistency cov}, we intend to prove the consistency in a similar way. However, there are two extra difficulties. One is that $\hat{y}_{ti} \not=\tilde{y}_{ti}$ and we need to control the error introduced by using $\hat{y}_{ti}$ in place of $\tilde{y}_{ti}$ in the estimators. The other is that the estimators $\check{F}_{t \mid k}(y)$ for $t=0,1$ and $k=1,\dots,K$ are not distribution functions and thus the representation theorem \citep{bobkov2019one}[Theorem \ref{thm:consistency cov}.11] can not be used directly. To solve this problem, we define
	\[
	\check{F}_{t\mid k}^*(y)=\sup_{v\leq y}\check{F}_{t \mid k}(v), 
	\]
	for $t=0,1$ and $k=1,\dots,K$. Then by defination, $\check{F}_{t\mid k}^*(y)$ is a distribution function and, for $u\in(0,1)$, $\check{F}_{t\mid k}^{*-1}(u)=
	\check{F}_{t\mid k}^{-1}(u)$.  Hence we can use $\check{F}_{t\mid k}^*(y)$ instead of $\check{F}_{t\mid k}(y)$ in the representation theorem. 
	
	We prove only for the lower bound, and the consistency result for the upper bound follows similarly. Let $\lambda_{k} = \pi_{\rm c}\pi_{k\mid \rm c}$. Because $\hat{\pi}_{\rm c}\hat{\pi}_{k\mid \rm c} = \hat{\lambda}_{1k}$, to prove $\check{\phi}^{2}_{\rm L} - \tilde{\phi}^{2}_{\rm L} \stackrel{P}{\to} 0$, we only need to prove
	\[
	\sum_{k=1}^{K}\hat{\lambda}_{1k}\int_0^{1}(\check{F}_{1\mid k}^{-1}(u)-\check{F}_{0\mid k}^{-1}(u))^2du - \sum_{k=1}^{K}\lambda_{k}\int_0^{1}(\tilde{F}_{1\mid k}^{-1}(u)-\tilde{F}_{0\mid k}^{-1}(u))^2du
	\stackrel{P}{\to} 0.
	\] 
	Note that
	\begin{align*}
		&\left\vert\sum_{k=1}^{K}\hat{\lambda}_{1k}\int_0^{1}(\check{F}_{1\mid k}^{-1}(u)-\check{F}_{0\mid k}^{-1}(u))^2du - \sum_{k=1}^{K}\lambda_{k}\int_0^{1}(\tilde{F}_{1\mid k}^{-1}(u)-\tilde{F}_{0\mid k}^{-1}(u))^2du\right\vert \\
		&\leq  \left\vert\sum_{k=1}^{K}\hat{\lambda}_{1k}\Big(\int_{0}^{1}(\check{F}_{1\mid k}^{-1}(u)-\check{F}_{0\mid k}^{-1}(u))^2du -  \int_{0}^{1}(\tilde{F}_{1\mid k}^{-1}(u) - \tilde{F}_{0\mid k}^{-1}(u))^2du\Big) \right\vert\\
		&+ \vert\sum_{k=1}^{K}(\hat{\lambda}_{1k}-\lambda_k)\int_{0}^{1}(\tilde{F}_{1\mid k}^{-1}(u) - \tilde{F}_{0\mid k}^{-1}(u))^2du\vert \\
		&\leq \left\vert\sum_{k=1}^K\Big(\frac{1}{n_1}\sum_{T_{i} = 1}d_{1i}1\{w_i = \xi_{k}\}\Big)\Big(\int_{0}^{1}(\check{F}_{1\mid k}^{-1}(u)-\check{F}_{0\mid k}^{-1}(u))^2du -  \int_{0}^{1}(\tilde{F}_{1\mid k}^{-1}(u) - \tilde{F}_{0\mid k}^{-1}(u))^2du\Big)\right\vert\\
		& + \left\vert\sum_{k=1}^K\Big(\frac{1}{n_0}\sum_{T_{i} = 0}d_{0i}1\{w_i = \xi_{k}\}\Big)\Big(\int_{0}^{1}(\check{F}_{1\mid k}^{-1}(u)-\check{F}_{0\mid k}^{-1}(u))^2du -  \int_{0}^{1}(\tilde{F}_{1\mid k}^{-1}(u) - \tilde{F}_{0\mid k}^{-1}(u))^2du\Big)\right\vert \\
		& + \left\vert\sum_{k=1}^{K}\left(\frac{\hat{\lambda}_{1k}}{\lambda_{k}}-1\right)\lambda_{1k}\int_{0}^{1}(\tilde{F}_{1\mid k}^{-1}(u) - \tilde{F}_{0\mid k}^{-1}(u))^2du\right\vert\\
		&\leq  2\max_k\left\vert\int_{0}^{1}(\check{F}_{1\mid k}^{-1}(u)-\check{F}_{0\mid k}^{-1}(u))^2du -  \int_{0}^{1}(\tilde{F}_{1\mid k}^{-1}(u) - \tilde{F}_{0\mid k}^{-1}(u))^2du\right\vert\\
		& + \max_k\left\vert\frac{\hat{\lambda}_{1k}}{\lambda_{k}} - 1\right\vert\sum_{k=1}^{K}\lambda_{1k}\int_{0}^{1}(\tilde{F}_{1\mid k}^{-1}(u) - \tilde{F}_{0\mid k}^{-1}(u))^2du\\
		&\equalscolon  I_{1,\rm c} + I_{2, \rm c}.
	\end{align*}
	
	By Condition \ref{cond:fourth moment} and Assumption 
	\ref{ass:iv}(iii), we have
	\[
	\theta \leq \frac{1}{N}\sum_{i=1}^{N}\vert y_{1i}\vert + \frac{1}{N}\sum_{i=1}^{N}\vert y_{0i}\vert \leq 2C_M^{\frac{1}{4}},
	\]
	and
	\[
	\vert\theta_{\rm c}\vert \leq C_{0}^{-1}\left(\frac{1}{N}\sum_{i=1}^{N}\vert y_{1i}\vert + \frac{1}{N}\sum_{i=1}^{N}\vert y_{0i}\vert\right)\leq 2C_{0}^{-1}C_M^{\frac{1}{4}}
	\]
	due to Jensen's inequality. Note that $\vert\tilde{y}_{ti}\vert \leq y_{ti} + \vert\theta_{\rm c}\vert$. By Condition \ref{cond:fourth moment} we have
	\begin{align*}
		\sum_{k=1}^{K}\lambda_{1k}\int_{0}^{1}(\tilde{F}_{1\mid k}^{-1}(u) - \tilde{F}_{0\mid k}^{-1}(u))^2du 
		&\leq \frac{2}{N}\sum_{i=1}^{N}\tilde{y}_{1i}^{2} + \frac{2}{N}\sum_{i=1}^{N}\tilde{y}_{0i}^{2} \\
		&\leq \frac{4}{N}\sum_{i=1}^{N}y_{1i}^{2} + \frac{4}{N}\sum_{i=1}^{N}y_{0i}^{2} + 8\vert\theta_{\rm c}\vert\\
		&\leq 8\sqrt{C_{M}}(1 + 2C_{0}^{-1}).
	\end{align*}
	For any small positive $\epsilon$, according to Hoeffding inequality for sampling without replacement \citep{bardenet2015} and the Bonferroni inequality we have
	\begin{align*}
		&\mathbb{P}\left( 8\sqrt{C_{M}}(1 + 2C_{0}^{-1})\max_{t,k}\frac{1}{\lambda_k}\vert\hat{\lambda}_{tk}-\lambda_k\vert\geq \epsilon\right)\\ &\leq 2K\exp\left( - \frac{1}{32C_{M}(1 + 2C_{0}^{-1})^2} n_1\land n_0\epsilon^2\min_k\lambda_k^2\right) \to 0
	\end{align*}
	by Conditions \ref{cond:positivity-compliance} and \ref{cond:order2}. Without loss of generality, we assume $\epsilon \leq 1/2$ in the proof.
	Thus for any $\delta > 0$ and sufficiently large $N$
	\[I_{2,\rm c} \leq \epsilon\]
	with probability at least $1-\delta/3$.
	
	Define $B_N = (C_{N} + C_{B})\lor 1$ where $C_{B} = 10C_{0}^{-2} C_M^{1/4}$.
	For $\epsilon \leq B_{N}(C_0\land C_{M}^{1/4})$, it is easy to verify that 
	\[
	\{B_N\vert\hat{\pi}_{\rm c} - \pi_{\rm c}\vert < \frac{\epsilon}{2}\}\cap\{B_N\vert\hat{\theta} - \theta\vert < \frac{\epsilon}{2}\} \subset \{B_N\vert\hat{\theta}_{\rm c} - \theta_{\rm c}\vert < \epsilon\}.
	\]
	Thus
	\begin{equation} \label{eq:prob bound theta_c}
		\begin{aligned}
			\mathbb{P}\left(B_N\vert\hat{\theta}_{\rm c} - \theta_{\rm c}\vert \geq C_{B}\epsilon \right) &\leq  \mathbb{P}\left(B_N\vert\hat{\pi}_{\rm c} - \pi_{\rm c}\vert \geq \frac{\epsilon}{2}\right) + \mathbb{P}\left(B_N\vert\hat{\theta} - \theta\vert\geq \frac{\epsilon}{2}\right)\\
			&\leq  4\exp\left(-\frac{2n_1\land n_0 \epsilon^2}{B_N^2}\right) + 
			4 \exp\left(-\frac{n_1\land n_0 \epsilon^2}{2B_N^2C_N^2}\right)
		\end{aligned}
	\end{equation}
	where the last inequality follows from the Hoeffding inequality. By Condition \ref{cond:order2}, 
	\[
	4\exp\left(-\frac{2n_1\land n_0 \epsilon^2}{B_N^2}\right) + 
	4 \exp\left(-\frac{n_1\land n_0 \epsilon^2}{2B_N^2C_N^2}\right) \to 0.
	\]
	Thus for sufficiently large $N$ with probability greater than $1 - \delta / 3$ we have $B_N\vert\hat{\theta}_{\rm c} - \theta_{\rm c}\vert \leq C_{B}\epsilon$.
	Note that $C_{0} \leq 1$, $\vert\tilde{y}_{ti}\vert \leq y_{ti} + \vert\theta_{\rm c}\vert \leq C_{N} + 2C_{0}^{-1}C_M^{1/4} \leq B_N$ and $\vert\hat{y}_{ti}\vert \leq \vert y_{ti}\vert + \vert\hat{\theta}_{\rm c}\vert \leq C_{N} + 2C_{0}^{-1}C_M^{1/4} + C_{B}\epsilon \leq B_N$ when the event $\{B_N\vert\hat{\theta}_{\rm c} - \theta_{\rm c}\vert < C_{B}\epsilon\}$ holds. 		
	On the event $\{B_N\vert\hat{\theta}_{\rm c} - \theta_{\rm c}\vert < C_{B}\epsilon\}$, by the representation theorem \citep{bobkov2019one}[Theorem \ref{thm:consistency cov}.11], using the similar arguments as those in the proof of Theorem \ref{thm:consistency cov}, we can show that
	\begin{align}
		&\max_k\vert\int_{0}^{1}(\check{F}_{1\mid k}^{-1}(u)-\check{F}_{0\mid k}^{-1}(u))^2du -  \int_{0}^{1}(\tilde{F}_{1\mid k}^{-1}(u) - \tilde{F}_{0\mid k}^{-1}(u))^2du \vert \notag\\
		&\leq 2\max_k\iint_{-B_N \leq v \leq w \leq B_N}[\vert\check{F}_{1\mid k}^{*}(v) - \tilde{F}_{1\mid k}(v)\vert + \vert\check{F}_{0\mid k}^{*}(v) - \tilde{F}_{0\mid k}(v)\vert \notag	\\ 
		& \phantom{\leq2\max_k\iint_{-B_N \leq v \leq w \leq B_N}} + \vert\check{F}_{1\mid k}^{*}(w) - \tilde{F}_{1\mid k}(w)\vert + \vert\check{F}_{0\mid k}^{*}(w) - \tilde{F}_{0\mid k}(w)\vert] dvdw \notag\\
		&\leq 4B_N\max_{k}\int_{-B_N}^{B_N}\vert\check{F}_{1\mid k}^{*}(v) - \tilde{F}_{1\mid k}(v)\vert dv + 4B_N\max_{k}\int_{-B_N}^{B_N}\vert\check{F}_{0\mid k}^{*}(v) - \tilde{F}_{0\mid k}(v)\vert dv \label{eq:upper bound I1c}.
	\end{align}	
	Here we only analyze the first term $B_N\max_{k}\int_{-B_N}^{B_N}\vert\check{F}_{1\mid k}^{*}(v) - F_{1\mid k}(v)\vert dv$ and the same result can be proved similarly for the second term.
	Define
	\[
	\bar{F}_{11\mid k}(y) = \frac{\hat{\lambda}_{1k}^{-1}}{n_1}\sum_{T_i=1}d_{1i}1\{\tilde{y}_{1i}\leq y\}1\{w_{i}=\xi_{k}\}, 
	\]
	\[
	\bar{F}_{01 \mid k}(y) = \frac{\hat{\lambda}_{1k}^{-1}}{n_0}\sum_{T_i=0}d_{0i}1\{\tilde{y}_{0i}\leq y\}1\{w_{i}=\xi_{k}\},
	\]
	\[
	\tilde{F}_{11\mid k}(y) = \frac{\lambda_{k}^{-1}}{N}\sum_{i=1}^Nd_{1i}1\{\tilde{y}_{1i}\leq y\}1\{w_{i}=\xi_{k}\}, 
	\]
	\[
	\tilde{F}_{01 \mid k}(y) = \frac{\lambda_{k}^{-1}}{N}\sum_{i=1}^Nd_{0i}1\{\tilde{y}_{0i}\leq y\}1\{w_{i}=\xi_{k}\}.
	\]
	for $k=1,\dots,K$. Let $\bar{F}_{1\mid k}(y) = \bar{F}_{11\mid k}(y) - \bar{F}_{01\mid k}(y)$, then the relationship 
	\[\check{F}_{1\mid k}(y) = \bar{F}_{1\mid k}(y + \hat{\theta}_{\rm c} - \theta_{\rm c})\] holds.
	Because $\vert\sup_{v\leq y} \bar{F}_{1\mid k}(v) - \sup_{v\leq y} \tilde{F}_{1\mid k}(v)\vert \leq \sup_{v}\vert\bar{F}_{1\mid k}(v) - \tilde{F}_{1\mid k}(v)\vert$ and $\sup_{v\leq y} \tilde{F}_{1\mid k}(v) = \tilde{F}_{1\mid k}(y)$, we have
	\begin{align*}
		&\vert\sup_{v \leq y}\bar{F}_{1\mid k}(v) - \bar{F}_{1\mid k}(y)\vert \\
		&\leq 
		\vert\sup_{v\leq y} \bar{F}_{1\mid k}(v) - \sup_{v\leq y} \tilde{F}_{1\mid k}(v)\vert + \vert \bar{F}_{1\mid k}(y) - \tilde{F}_{1\mid k}(y)\vert\\
		&\leq 2\sup_{v}\vert\bar{F}_{1\mid k}(v) - \tilde{F}_{1\mid k}(v)\vert.
	\end{align*}
	Hence
	\begin{align*}
		&\vert\check{F}_{1\mid k}^{*}(y) - \tilde{F}_{1\mid k}(y)\vert \\
		& \leq 
		\vert\sup_{v \leq y + \hat{\theta}_{\rm c} - \theta_{\rm c}}\bar{F}_{1\mid k}(v) - \bar{F}_{1\mid k}(y + \hat{\theta}_{\rm c} - \theta_{\rm c})\vert  + \vert\bar{F}_{1\mid k}(y + \hat{\theta}_{\rm c} - \theta_{\rm c}) - \tilde{F}_{1\mid k}(y)\vert \\
		&\leq 2\sup_{v}\vert\bar{F}_{1\mid k}(v) - \tilde{F}_{1\mid k}(v)\vert
		+ \vert\bar{F}_{1\mid k}(y + \hat{\theta}_{\rm c} - \theta_{\rm c}) - \tilde{F}_{1\mid k}(y + \hat{\theta}_{\rm c} - \theta_{\rm c})\vert \\
		& + \vert\tilde{F}_{1\mid k}(y + \hat{\theta}_{\rm c} - \theta_{\rm c}) - \tilde{F}_{1\mid k}(y)\vert \\
		&\leq 3\sup_{v}\vert\bar{F}_{1\mid k}(v) - \tilde{F}_{1\mid k}(v)\vert
		+ \vert\tilde{F}_{1\mid k}(y + \hat{\theta}_{\rm c} - \theta_{\rm c}) - \tilde{F}_{1\mid k}(y)\vert.
	\end{align*}
	Moreover, because 
	\begin{align*}
		&\vert\tilde{F}_{1\mid k}(y + \hat{\theta}_{\rm c} - \theta_{\rm c}) - \tilde{F}_{1\mid k}(y)\vert\\ 
		&\leq  \vert\tilde{F}_{11\mid k}(y + \hat{\theta}_{\rm c} - \theta_{\rm c}) -  \tilde{F}_{11\mid k}(y)\vert + \vert\tilde{F}_{01\mid k}(y + \hat{\theta}_{\rm c} - \theta_{\rm c}) -  \tilde{F}_{01\mid k}(y)\vert \\
		&\leq  \lambda^{-1}_{k}\Big(\frac{1}{N}\sum_{i=1}^{N}d_{1i}1\{\vert\tilde{y}_{1i} - y\vert\leq \vert\hat{\theta}_{\rm c} - \theta_{\rm c}\vert\}1\{ w_{i} = \xi_{k}\}\\
		& \phantom{\leq  \lambda^{-1}_{k}\Big(} + \frac{1}{N}\sum_{i=1}^{N}d_{0i}1\{\vert\tilde{y}_{0i} - y\vert\leq \vert\hat{\theta}_{\rm c} - \theta_{\rm c}\vert\}1\{ w_{i} = \xi_{k}\}\Big),
	\end{align*}
	we have
	\[
	B_N\max_{k}\int_{-B_N}^{B_N}\vert\check{F}_{1\mid k}^{*}(v) - F_{1\mid k}(v)\vert dv \leq 3B_N^2 \max_{k}\sup_{v}\vert\bar{F}_{1\mid k}(v) - \tilde{F}_{1\mid k}(v)\vert + 2B_N\vert\hat{\theta}_{\rm c} - \theta_{\rm c}\vert.
	\]
	By inequality \eqref{eq:prob bound theta_c}, we have $B_N\vert\hat{\theta}_{\rm c} - \theta_{\rm c}\vert \leq C_{B}\epsilon$ with probability at least $1 - \delta / 3$ for sufficiently large $N$.
	Using the similar arguments as those used to analyze $I_{31}$ in the proof of Theorem \ref{thm:consistency cov}, we can show
	\[B_N^2 \max_{k}\sup_{v}\vert\bar{F}_{1\mid k}(v) - \tilde{F}_{1\mid k}(v)\vert \leq \epsilon\]
	with probability at least $1 - \delta/6$ for sufficiently large $N$. By applying the similar arguments to the second term of expression \eqref{eq:upper bound I1c}, we have $I_{1,\rm c} \leq (48+ 32C_B)\epsilon$ with probability at least $1 - 2\delta/3$ for sufficiently large $N$. Thus we have proved that for any small positive numbers $\epsilon$ and $\delta$, we have $\vert\check{\phi}^{2}_{\rm L} - \tilde{\phi}^{2}_{\rm L}\vert \leq (49 + 32C_{B})\epsilon$ with probability at least $1 - \delta$ for sufficiently large $N$ and this implies the consistency of the estimator.
\end{proof}

\section{Further simulation results}
To explore the reliability of the proposed confidence intervals (CIs), we consider the case where $\phi^{2}(\tau)$ attains the lower bound $\phi_{\rm L}^{2}$. In this case, the asymptotic variance of $\sqrt{N}(\hat{\theta} - \theta)$ achieves its upper bound. According to the proof of the sharpness in Theorem \ref{thm:cov}, we can modify the finite populations considered in Section \ref{subsec: sim cov} to make $\phi^{2}(\tau)$ equal to $\phi_{\rm L}^{2}$ without changing $\pi_{k}$, $F_{1\mid k}(y)$ or $F_{0\mid k}(y)$ ($k=1,\dots,K$). Then we conducted the simulation in the same way as in Section \ref{subsec: sim cov} in the main text under the modified finite populations. The average width (AW) and coverage rate (CR) of $95\%$ CIs based on the naive lower bound zero \citep{neyman1990application}, the estimator of $\phi^2_{\rm AL}$ \citep{Aronow2014SBotViRE}, the estimator of $\phi^2_{\rm DL}$ \citep{Ding2018DTEV} and the estimator of $\phi^2_{\rm L}$ are summarized in the following table.

\begin{table}[h]
	\def~{\hphantom{0}}
	\centering
	\caption{Average widths (AWs) and coverage rates (CRs) of 95\% CIs based on the naive bound, $\phi^2_{\rm AL}$, $\phi^2_{\rm DL}$ and $\phi^2_{\rm L}$ under different population sizes when $\phi^{2}(\tau)$ attains the lower bound $\phi_{\rm L}^{2}$ ($n_{1} = n_{0} = N/2$)}
	\vspace{10pt}
		\begin{tabular}{*{9}{c}}
			\toprule
			\multirow{2}{*}{Method}&\multicolumn{2}{c}{naive}&\multicolumn{2}{c}{$\phi^2_{\rm AL}$}&\multicolumn{2}{c}{$\phi^2_{\rm DL}$}
			&\multicolumn{2}{c}{$\phi^2_{\rm L}$}\\
			&{\footnotesize AW} & {\footnotesize CR} &{\footnotesize AW} & {\footnotesize CR} &{\footnotesize AW} & {\footnotesize CR} &{\footnotesize AW} &{\footnotesize CR}\\
			\midrule
			$N = 400$ & 1.511 & 96.4\% & 1.495 & 96.3\% & 1.493 & 96.3\% & 1.380 & 94.7\%\\
			$N = 800$ & 1.033 & 96.6\% & 1.025 & 96.4\% & 1.025 & 96.4\% & 0.944 & 95.2\%\\
			$N = 2000$ & 0.674 & 97.0\% & 0.669 & 96.9\% & 0.669 & 96.9\% & 0.619 & 95.7\%\\
			\bottomrule
		\end{tabular}
	\label{table: CI cov Ated}
\end{table}

Comparing Table \ref{table: CI cov Ated} with Table \ref{table: CI cov} in the main text, we find that the AWs are similar while the CRs are smaller in Table \ref{table: CI cov Ated}. This is because the bounds are all the same under the finite populations considered in Table \ref{table: CI cov Ated} and Table \ref{table: CI cov}, but the variance of $\hat{\theta}$ is larger here. The CI based on $\phi_{\rm L}^{2}$ is still quite reliable in this case because its CR is close to $95\%$. 

We then investigate the AW and CR of the CIs in randomized experiments with noncompliance similarly. The following table summarizes the average width (AW) and coverage rate (CR) of $95\%$ CIs based on the naive lower bound zero, and the estimator of lower bounds in Theorem \ref{thm:compliance} using and without using the covariate. 

\begin{table}[h]
	\def~{\hphantom{0}}
	\centering
	\caption{Average widths (AWs) and coverage rates (CRs) of 95\% CIs based on the naive bound, HNL and HL under different population sizes when $\phi^{2}(\tilde{\tau})$ attains the lower bound $\tilde{\phi}_{\rm L}^{2}$.  LNC: lower bound without covariate; HNC: upper bound without covariate; LC: lower bound with covariate; HC: upper bound with covariate ($n_{1} = n_{0} = N/2$)}
	\vspace{10pt}
		\begin{tabular}{*{7}{c}}
			\toprule
			\multirow{2}{*}{Method}&\multicolumn{2}{c}{naive}&\multicolumn{2}{c}{LNC}&\multicolumn{2}{c}{LC}\\
			&{\footnotesize AW} & {\footnotesize CR} &{\footnotesize AW} & {\footnotesize CR} &{\footnotesize AW} & {\footnotesize CR}\\
			\midrule
			$N = 400$ & 2.186 & 96.7\% & 2.166 & 96.7\% & 1.983 & 94.1\%\\
			$N = 800$ & 1.553 & 97.0\% & 1.542 & 96.8\% & 1.397 & 94.5\%\\
			$N = 2000$ & 0.980 & 96.2\% & 0.973 & 96.0\% & 0.893 & 94.8\%\\
			\bottomrule
		\end{tabular}
	\label{table: CI compliance Ated}
\end{table}

It can be seen that the CI based on LC is still quite reliable even when the asymptotic variance of $\hat{\theta}_{\rm c}$ achieves its upper bound.

\bibliographystyle{chicago}
\bibliography{vb}
\end{document}